\documentclass[11pt,a4paper,reqno]{amsart}

\usepackage{color,latexsym,amsfonts,amssymb,amsmath,cite, amsthm, bbm}
\usepackage{hyperref}
\usepackage{fullpage}
\usepackage{fancyhdr}
\usepackage{graphicx}
\usepackage{tikz}
\usepackage{pgfplots}

\definecolor{wiasblue}   {cmyk}{1.0, 0.60, 0, 0}

\def\Z{\mathbb{Z}}
\def\E{\mathbb{E}}
\def\P{\mathbb{P}}

\def\R{\mathbb{R}}
\def\S{\sigma}
\def\bR{\breve{\mathbb{R}}}
\def\mc{\mathcal}
\def\ms{\mathsf}

\def\la{\lambda}
\def\g{\gamma}
\def\de{\delta}
\def\es{\emptyset}
\def\one{\mathbbmss{1}}

\def\d{{\rm d}}
\def\b{\beta}
\def\X{\mc X}

\def\e{\varepsilon}
\def\V{V}

\def\CO{\mc C}

\def\PD{\beta}
\def\PDD{\ms{PD}^{M, q}}

\def\bim{B^M_i}
\def\dim{D^M_i}

\def\PDDn{\ms{PD}^{M, q}(\PP_n)}

\def\PDodn{\b^{M, 0}_d(\PP_n)}
\def\PDoodn{\b^{M, 1}_{b, d}(\PP_n)}
\def\mcp{\ms{MatC}}
\def\str{\ms{Str}}

\def\x{x}
\def\Poi{\ms{Poi}}
\def\kk{\{1, \dots, k\}}
\def\four{\{1, 2, 3, 4\}}
\def\ckn{c^k_n}
\def\cfn{c^4_n}
\def\bx{\breve{x}}
\def\xx{\boldsymbol x}
\def\yy{\boldsymbol y}
\def\zz{\boldsymbol z}

\def\bxx{\breve{\xx}}
\def\ff{\boldsymbol f}
\def\mm{E''}

\def\APFz{\ms{APF}^{M, 0}}
\def\APFo{\ms{APF}^{M, 1}}

\def\mPD{\overline{\beta}}
\def\Var{\ms{Var}}
\def\Wo{W_1}
\def\Wn{W_n}
\def\Wmn{W^{(1,2)}_n}
\def\Wr{W_r}

\def\bWn{\breve{W}_n}

\def\mun{\mu_n}
\def\nun{\nu_n}
\def\tf{r_{\ms f}}
\def\tfo{[0, r_{\ms f}]}
\def\tfs{[0, r_{\ms f}]^2}

\def\PP{\mc P}
\def\PPn{\mc P_n}

\def\bPPn{\breve{\mc P_n}}

\def\lan{\langle}
\def\ran{\rangle}

\def\dist{\ms{dist}}
\def\rac{r_{\ms{AC}}}

\def\Tc{T_{\ms C}}
\def\Tl{T_{\ms L}}
\def\rc{r_{\ms C}}
\def\rl{r_{\ms L}}
\def\La{\Lambda}
\def\NN{\mc N}

\newtheorem{theorem}{Theorem}[section]
\newtheorem{corollary}[theorem]{Corollary}

\newtheorem{lemma}[theorem]{Lemma}
\newtheorem{proposition}[theorem]{Proposition}
\theoremstyle{definition}

\keywords{Point processes, goodness-of-fit tests, central limit theorem, topological data analysis, persistent Betti number}
\subjclass[2010]{60D05; 55N20; 60F17}

\begin{document}

\author{Christophe A.N.~Biscio}
\author{Nicolas Chenavier}
\author{Christian Hirsch}
\author{Anne Marie Svane}
\address[Christophe A.N.~Biscio, Anne Marie Svane]{Aalborg University, Department of Mathematical Sciences,
         Skjernvej 4,  9220 Aalborg \O, Denmark}
	 \email{christophe@math.aau.dk, annemarie@math.aau.dk}
\address[Nicolas Chenavier]{Universit\'e Littoral C{\^o}te d'Opale, 
         EA 2797, LMPA, 50 rue Ferdinand Buisson, 62228 Calais, France}
\email{nicolas.chenavier@univ-littoral.fr}
\address[Christian Hirsch]{University of Mannheim, Institute of Mathematics, 68161 Mannheim, Germany}
\email{hirsch@uni-mannheim.de}

	\title{Testing goodness of fit for point processes via topological data analysis}

\date{\today}
\begin{abstract}
        We introduce tests for the goodness of fit of point patterns via methods from topological data analysis. More precisely, the persistent Betti numbers give rise to a bivariate functional summary statistic for observed point patterns that is asymptotically Gaussian in large observation windows.  We analyze the power of tests derived from this statistic on simulated point patterns and compare its performance with global envelope tests. Finally, we apply the tests to a point pattern from an application context in neuroscience. As the main methodological contribution, we derive sufficient conditions for a functional central limit theorem on bounded persistent Betti numbers of point processes with exponential decay of correlations.
\end{abstract}

\maketitle

\section{Introduction}
\label{intrSec}
\begin{figure}[!htpb]
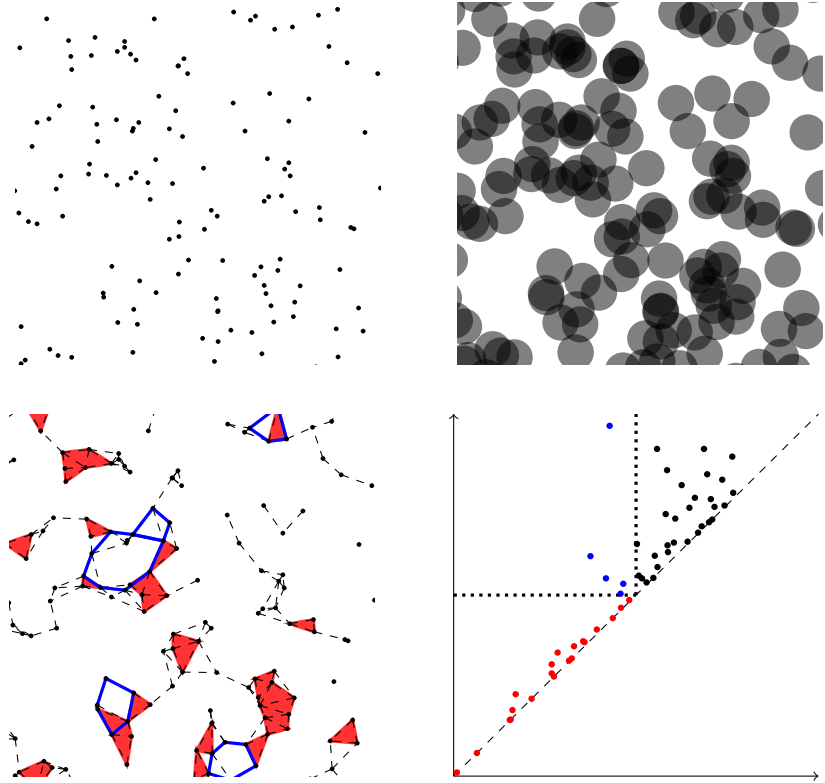

     \begin{tikzpicture}[scale=.60]
 \clip(-4, -4) rectangle (4, 4);
\fill[black] (2.09, -4.36) circle (1.5pt);\fill[black] (-0.62, 1.21) circle (1.5pt);\fill[black] (-3.00, 0.17) circle (1.5pt);\fill[black] (2.67, -0.53) circle (1.5pt);\fill[black] (0.13, 0.34) circle (1.5pt);\fill[black] (-4.55, -3.09) circle (1.5pt);\fill[black] (2.00, 1.38) circle (1.5pt);\fill[black] (1.46, -2.25) circle (1.5pt);\fill[black] (-4.58, 3.65) circle (1.5pt);\fill[black] (-3.92, -0.65) circle (1.5pt);\fill[black] (-1.09, 0.02) circle (1.5pt);\fill[black] (-2.30, 1.69) circle (1.5pt);\fill[black] (1.41, -3.76) circle (1.5pt);\fill[black] (-4.22, -4.85) circle (1.5pt);\fill[black] (-2.23, 4.20) circle (1.5pt);\fill[black] (1.76, -1.92) circle (1.5pt);\fill[black] (3.35, -4.16) circle (1.5pt);\fill[black] (-1.35, 0.24) circle (1.5pt);\fill[black] (-4.26, -4.95) circle (1.5pt);\fill[black] (-3.31, 3.64) circle (1.5pt);\fill[black] (1.16, 3.61) circle (1.5pt);\fill[black] (-3.06, -3.80) circle (1.5pt);\fill[black] (2.10, -2.94) circle (1.5pt);\fill[black] (3.43, -1.00) circle (1.5pt);\fill[black] (-1.91, -4.70) circle (1.5pt);\fill[black] (1.50, -2.43) circle (1.5pt);\fill[black] (0.78, 1.78) circle (1.5pt);\fill[black] (-2.43, -4.41) circle (1.5pt);\fill[black] (-1.76, -0.03) circle (1.5pt);\fill[black] (-1.32, -3.11) circle (1.5pt);\fill[black] (-1.42, -3.73) circle (1.5pt);\fill[black] (-0.91, -2.41) circle (1.5pt);\fill[black] (-0.62, -1.23) circle (1.5pt);\fill[black] (1.19, -3.29) circle (1.5pt);\fill[black] (-2.78, 2.81) circle (1.5pt);\fill[black] (1.27, 3.70) circle (1.5pt);\fill[black] (-3.78, -3.90) circle (1.5pt);\fill[black] (2.02, 4.58) circle (1.5pt);\fill[black] (-1.27, -2.15) circle (1.5pt);\fill[black] (-0.41, -1.17) circle (1.5pt);\fill[black] (-1.37, -4.66) circle (1.5pt);\fill[black] (3.36, -0.97) circle (1.5pt);\fill[black] (-2.21, 1.31) circle (1.5pt);\fill[black] (-4.07, -3.02) circle (1.5pt);\fill[black] (-4.44, -3.20) circle (1.5pt);\fill[black] (2.00, 0.46) circle (1.5pt);\fill[black] (4.54, -3.76) circle (1.5pt);\fill[black] (-0.42, 2.60) circle (1.5pt);\fill[black] (2.13, -2.62) circle (1.5pt);\fill[black] (3.04, -3.21) circle (1.5pt);\fill[black] (4.96, 2.09) circle (1.5pt);\fill[black] (2.65, 3.72) circle (1.5pt);\fill[black] (1.39, -1.84) circle (1.5pt);\fill[black] (-2.00, 0.18) circle (1.5pt);\fill[black] (-1.41, 1.22) circle (1.5pt);\fill[black] (-0.74, 4.30) circle (1.5pt);\fill[black] (2.15, -0.62) circle (1.5pt);\fill[black] (0.26, 4.48) circle (1.5pt);\fill[black] (0.11, -3.42) circle (1.5pt);\fill[black] (-2.33, 2.83) circle (1.5pt);\fill[black] (0.79, 3.79) circle (1.5pt);\fill[black] (1.56, -0.29) circle (1.5pt);\fill[black] (-3.03, 1.72) circle (1.5pt);\fill[black] (0.86, -4.08) circle (1.5pt);\fill[black] (-3.71, -0.84) circle (1.5pt);\fill[black] (-4.62, 3.61) circle (1.5pt);\fill[black] (2.25, -2.26) circle (1.5pt);\fill[black] (-4.38, 2.68) circle (1.5pt);\fill[black] (-3.63, 0.82) circle (1.5pt);\fill[black] (4.44, -1.19) circle (1.5pt);\fill[black] (0.08, -4.52) circle (1.5pt);\fill[black] (4.95, -1.53) circle (1.5pt);\fill[black] (1.61, -2.70) circle (1.5pt);\fill[black] (-4.26, -4.96) circle (1.5pt);\fill[black] (4.67, -2.20) circle (1.5pt);\fill[black] (-3.52, -0.89) circle (1.5pt);\fill[black] (1.58, 1.99) circle (1.5pt);\fill[black] (4.56, 4.75) circle (1.5pt);\fill[black] (1.53, 0.55) circle (1.5pt);\fill[black] (4.57, -2.20) circle (1.5pt);\fill[black] (0.44, -2.83) circle (1.5pt);\fill[black] (-4.83, 3.56) circle (1.5pt);\fill[black] (-1.60, 3.14) circle (1.5pt);\fill[black] (-0.30, 2.76) circle (1.5pt);\fill[black] (0.45, -2.81) circle (1.5pt);\fill[black] (-4.81, 2.22) circle (1.5pt);\fill[black] (-4.98, 0.38) circle (1.5pt);\fill[black] (1.81, -1.67) circle (1.5pt);\fill[black] (0.38, -4.68) circle (1.5pt);\fill[black] (1.46, -0.39) circle (1.5pt);\fill[black] (1.91, 0.32) circle (1.5pt);\fill[black] (1.89, 4.17) circle (1.5pt);\fill[black] (3.29, 4.86) circle (1.5pt);\fill[black] (2.60, -4.05) circle (1.5pt);\fill[black] (-4.51, 2.06) circle (1.5pt);\fill[black] (-2.06, -2.51) circle (1.5pt);\fill[black] (-4.31, -2.53) circle (1.5pt);\fill[black] (3.07, 3.85) circle (1.5pt);\fill[black] (-1.88, -1.83) circle (1.5pt);\fill[black] (1.69, 3.41) circle (1.5pt);\fill[black] (0.06, -4.06) circle (1.5pt);\fill[black] (4.65, -2.12) circle (1.5pt);\fill[black] (1.32, -4.57) circle (1.5pt);\fill[black] (-4.92, 0.76) circle (1.5pt);\fill[black] (-4.75, 1.87) circle (1.5pt);\fill[black] (3.63, -3.26) circle (1.5pt);\fill[black] (-4.90, -2.86) circle (1.5pt);\fill[black] (1.15, 4.57) circle (1.5pt);\fill[black] (3.40, 4.56) circle (1.5pt);\fill[black] (-1.47, 2.86) circle (1.5pt);\fill[black] (3.78, 4.39) circle (1.5pt);\fill[black] (-4.28, -0.71) circle (1.5pt);\fill[black] (-3.07, 4.31) circle (1.5pt);\fill[black] (-3.87, -3.98) circle (1.5pt);\fill[black] (-4.18, -1.69) circle (1.5pt);\fill[black] (-1.78, 1.41) circle (1.5pt);\fill[black] (-4.09, 1.15) circle (1.5pt);\fill[black] (2.90, 4.88) circle (1.5pt);\fill[black] (-0.55, 0.00) circle (1.5pt);\fill[black] (4.86, -2.15) circle (1.5pt);\fill[black] (1.25, -4.72) circle (1.5pt);\fill[black] (-4.81, 3.81) circle (1.5pt);\fill[black] (0.37, -3.93) circle (1.5pt);\fill[black] (0.55, -1.70) circle (1.5pt);\fill[black] (0.15, -1.01) circle (1.5pt);\fill[black] (-0.40, -0.85) circle (1.5pt);\fill[black] (-2.37, 0.43) circle (1.5pt);\fill[black] (-4.57, 0.71) circle (1.5pt);\fill[black] (-4.24, 3.82) circle (1.5pt);\fill[black] (-2.43, 4.98) circle (1.5pt);\fill[black] (-1.40, -2.78) circle (1.5pt);\fill[black] (-3.10, -0.12) circle (1.5pt);\fill[black] (3.08, -3.83) circle (1.5pt);\fill[black] (-0.42, -2.00) circle (1.5pt);\fill[black] (-4.01, -0.16) circle (1.5pt);\fill[black] (-2.25, 3.23) circle (1.5pt);\fill[black] (0.43, -0.72) circle (1.5pt);\fill[black] (-1.76, 2.92) circle (1.5pt);\fill[black] (-3.24, -3.66) circle (1.5pt);\fill[black] (-0.23, -4.49) circle (1.5pt);\fill[black] (0.41, -2.54) circle (1.5pt);\fill[black] (-3.39, -4.28) circle (1.5pt);\fill[black] (-2.39, 0.21) circle (1.5pt);\fill[black] (-3.26, 1.43) circle (1.5pt);\fill[black] (-4.52, 2.79) circle (1.5pt);\fill[black] (-0.85, 1.93) circle (1.5pt);\fill[black] (2.08, 3.46) circle (1.5pt);\fill[black] (1.25, -2.03) circle (1.5pt);\fill[black] (0.30, -0.59) circle (1.5pt);\fill[black] (-0.92, 0.53) circle (1.5pt);\fill[black] (-2.83, 3.16) circle (1.5pt);\fill[black] (-2.76, 2.52) circle (1.5pt);\fill[black] (-1.04, -0.35) circle (1.5pt);\fill[black] (4.78, -4.82) circle (1.5pt);\fill[black] (3.28, 2.67) circle (1.5pt);\fill[black] (-0.38, -4.13) circle (1.5pt);\fill[black] (-3.52, 1.35) circle (1.5pt);\fill[black] (-2.07, -2.42) circle (1.5pt);\fill[black] (1.60, 4.93) circle (1.5pt);\fill[black] (-4.61, -1.32) circle (1.5pt);\fill[black] (-3.63, 4.19) circle (1.5pt);\fill[black] (-2.76, -3.84) circle (1.5pt);\fill[black] (3.53, -2.75) circle (1.5pt);\fill[black] (2.44, 1.87) circle (1.5pt);\fill[black] (1.84, -0.27) circle (1.5pt);\fill[black] (0.45, -3.49) circle (1.5pt);\fill[black] (-1.46, 0.12) circle (1.5pt);\fill[black] (-1.29, 1.35) circle (1.5pt);\fill[black] (3.67, 1.13) circle (1.5pt);\fill[black] (-0.48, 1.61) circle (1.5pt);\fill[black] (3.12, -1.90) circle (1.5pt);\fill[black] (2.69, -0.80) circle (1.5pt);\fill[black] (-1.60, 3.01) circle (1.5pt);\fill[black] (-2.96, -0.72) circle (1.5pt);\fill[black] (-0.22, 2.43) circle (1.5pt);\fill[black] (-4.50, -2.17) circle (1.5pt);\fill[black] (-4.79, 4.65) circle (1.5pt);\fill[black] (3.81, -4.55) circle (1.5pt);\fill[black] (-0.20, -1.69) circle (1.5pt);\fill[black] (4.96, -4.81) circle (1.5pt);\fill[black] (3.95, 2.43) circle (1.5pt);\fill[black] (-1.43, 1.16) circle (1.5pt);\fill[black] (4.01, -0.10) circle (1.5pt);\fill[black] (1.04, 1.16) circle (1.5pt);\fill[black] (1.67, -3.06) circle (1.5pt);\fill[black] (2.88, 3.03) circle (1.5pt);\fill[black] (4.72, 4.48) circle (1.5pt);\fill[black] (-0.43, 2.60) circle (1.5pt);\fill[black] (-0.33, 1.04) circle (1.5pt);\fill[black] (-1.75, -3.07) circle (1.5pt);\fill[black] (4.29, -0.90) circle (1.5pt);\fill[black] (-3.90, 3.00) circle (1.5pt);\fill[black] (-1.35, 2.73) circle (1.5pt);\fill[black] (4.88, -3.67) circle (1.5pt);\fill[black] (0.73, -3.24) circle (1.5pt);\fill[black] (-2.19, 0.93) circle (1.5pt);\fill[black] (-3.87, -3.16) circle (1.5pt);\fill[black] (2.03, 0.17) circle (1.5pt);\fill[black] (4.64, -2.12) circle (1.5pt);\fill[black] (-0.95, 3.64) circle (1.5pt);
\end{tikzpicture} \qquad 
     \input{poisball}  \\ \vspace{0.6cm} 	
     \input{poisalphagraph} \qquad  
     \begin{tikzpicture}[scale=.60]
\draw[->] (0, 0) -- (
8
, 0);
\draw[->] (0, 0) -- (0, 
8
);
\draw[dashed] (0, 0) -- (
8
, 
8
);
\draw[dotted, very thick] (4,4) -- (4, 8);
\draw[dotted, very thick] (0, 4) -- (4, 4);
\fill[blue] (3.00, 4.86) circle (2.0pt);\fill[blue] (3.34, 4.37) circle (2.0pt);\fill[blue] (3.42, 7.74) circle (2.0pt);\fill[blue] (3.66, 4.03) circle (2.0pt);\fill[blue] (3.72, 4.25) circle (2.0pt);
\fill[red] (3.49, 3.49) circle (2.0pt);\fill[red] (1.25, 1.25) circle (2.0pt);\fill[red] (1.29, 1.46) circle (2.0pt);\fill[red] (2.53, 2.55) circle (2.0pt);\fill[red] (1.71, 1.71) circle (2.0pt);\fill[red] (0.07, 0.08) circle (2.0pt);\fill[red] (2.64, 2.87) circle (2.0pt);\fill[red] (2.84, 2.98) circle (2.0pt);\fill[red] (3.14, 3.24) circle (2.0pt);\fill[red] (2.15, 2.47) circle (2.0pt);\fill[red] (2.20, 2.20) circle (2.0pt);\fill[red] (3.67, 3.72) circle (2.0pt);\fill[red] (2.28, 2.73) circle (2.0pt);\fill[red] (2.15, 2.27) circle (2.0pt);\fill[red] (1.23, 1.24) circle (2.0pt);\fill[red] (1.36, 1.81) circle (2.0pt);\fill[red] (2.87, 2.96) circle (2.0pt);\fill[red] (2.59, 2.60) circle (2.0pt);\fill[red] (3.85, 3.89) circle (2.0pt);\fill[red] (0.51, 0.51) circle (2.0pt);
\fill[black] (4.47, 4.62) circle (2.0pt);\fill[black] (5.94, 5.98) circle (2.0pt);\fill[black] (6.11, 7.06) circle (2.0pt);\fill[black] (5.44, 5.52) circle (2.0pt);\fill[black] (5.60, 5.61) circle (2.0pt);\fill[black] (6.13, 6.26) circle (2.0pt);\fill[black] (4.13, 4.37) circle (2.0pt);\fill[black] (5.67, 5.67) circle (2.0pt);\fill[black] (4.46, 7.23) circle (2.0pt);\fill[black] (5.29, 6.15) circle (2.0pt);\fill[black] (5.13, 5.18) circle (2.0pt);\fill[black] (5.49, 7.23) circle (2.0pt);\fill[black] (5.89, 6.55) circle (2.0pt);\fill[black] (4.23, 4.28) circle (2.0pt);\fill[black] (4.68, 6.76) circle (2.0pt);\fill[black] (4.41, 4.87) circle (2.0pt);\fill[black] (5.17, 5.93) circle (2.0pt);\fill[black] (4.69, 5.10) circle (2.0pt);\fill[black] (5.00, 6.43) circle (2.0pt);\fill[black] (5.64, 6.12) circle (2.0pt);\fill[black] (5.56, 6.67) circle (2.0pt);\fill[black] (4.65, 5.79) circle (2.0pt);\fill[black] (4.38, 4.38) circle (2.0pt);\fill[black] (5.72, 5.95) circle (2.0pt);\fill[black] (4.02, 5.13) circle (2.0pt);\fill[black] (4.86, 5.69) circle (2.0pt);\fill[black] (5.36, 5.38) circle (2.0pt);\fill[black] (4.70, 4.95) circle (2.0pt);\fill[black] (4.06, 4.43) circle (2.0pt);\fill[black] (4.83, 5.16) circle (2.0pt);
\end{tikzpicture}
     \caption{Top: Realization of Poisson point process (left) and union of balls centered at the points of the process (right). Bottom: Alpha-complex corresponding to the union of balls with alive (blue) and dead (red) loops marked (left). Associated persistence diagram (right).}
     \label{persFig}
\end{figure}

Topological data analysis (TDA) provides insights into a variety of datasets by capturing their most salient properties via refined topological features.
Since the mathematical field of topology specializes in describing invariants of objects independently of the choice of a precise metric, these features are robust against small perturbations or different embeddings of the object \cite{carlsson, chazal}.
Among the most classical topological invariants are the Betti numbers.
Loosely speaking, they capture the number of $k$-dimensional holes of the investigated structure.
TDA refines this idea substantially by constructing filtrations and tracing when topological features appear and disappear.
In point pattern analysis, simplicial complexes are built so that they are topologically equivalent
to a union of balls with the same radius and centered at the data points,
see the first three panels of Figure~\ref{persFig}.
As the radius increases,
a sequence of simplicial complexes is then defined.
Examples of such complexes 
are the basic \v{C}ech complex or the more elaborate $\alpha$-complex,
which is based on the Delaunay triangulation, see \cite{edHar}.
In that framework, 1-dimensional
features correspond to loops in the simplicial complexes while 0-dimensional features correspond to connected components.
When moving up in the filtration, additional edges appear and at some point create new loops.
On the other hand, more and more triangles also appear, thereby causing completely filled loops to disappear.
Usually, the filtration is indexed by time, and we refer to the appearance and disappearance of features as births and deaths. We refer the reader to \cite{edHar} for a detailed presentation of these concepts.
The \emph{persistence diagram} visualizes the time points
when the features are born and die, see the bottom-right panel in Figure \ref{persFig}.
Persistent Betti numbers count the number of events in upper-left blocks of the persistence diagram and are also illustrated in the figure.

In this paper, we leverage persistent Betti numbers to derive goodness-of-fit tests for planar point processes. In this setting, the abstract general definition of persistent Betti numbers gives way to a clear geometric intuition induced by a picture of growing disks centered at the points of the pattern and all having radius $r$, corresponding to the index of the filtration. Features of dimension 0 correspond to connected components in the union of balls, interpreted as point clusters, whereas boundaries of the complement set can be considered as the loops forming the 1-dimensional features.
Since the notion of clusters in the sense of connected components lies at the heart of persistent Betti numbers in degree 0, they become highly attractive as a tool to detect clustering in point patterns. Our tests are based on a novel functional central limit theorem (CLT) for the persistent Betti numbers in large domains. The present work embeds into two active streams of current research.

%
%
First, now that TDA has become widely adopted, the community is vigorously working towards putting the approach on a firm statistical foundation paving the way for hypothesis testing. On the one hand, this encompasses large-sample Monte Carlo tests when working on a fixed domain \cite{biscioTDA, bubenik, silhouette}. Although these tests are highly flexible, the test statistics under the null hypothesis must be re-computed each time when testing observations in a different window. In large domains, this becomes time-consuming. On the other hand, there has been substantial progress towards establishing CLTs in large domains for functionals related to persistent Betti numbers \cite{yogeshAdler1, yogeshAdler2, krebs, owada, shirai}. However, these results are restricted to the null hypothesis of complete spatial randomness -- i.e., the Poisson point process -- and establish asymptotic Gaussianity on a multivariate, but not on a functional level. Our proof of a functional CLT is based on recently developed stabilization techniques for point processes with exponential decay of correlations \cite{yogeshCLT}. As explained in the final section of \cite{calka}, the main technical step towards a functional CLT are bounds on the cumulants.

%
%
Second, the introduction of global rank envelope tests has lead to a novel surge of research activity in goodness-of-fit tests for point processes \cite{envelop}. One of the reasons for their popularity is that they rely on functional summary statistics rather than scalar quantities. Thus, they reveal a substantially more fine-grained picture of the underlying point pattern. In the overwhelming majority of cases, variants of the $K$-function are used as a functional summary statistic, thereby essentially capturing the relative density of point pairs at different distances. Here, the persistent Betti numbers offer an opportunity to augment the basic second-order information by more refined characteristics of the data. Still, even for classical summary statistics, rigorous limit theorems in large domains remain scarce. For instance, a functional central limit theorem of the estimated $K$-function is proven in detail only for the Poisson point process in \cite{kclt} and an extension to $\alpha$-determinantal point processes is outlined in \cite{alphaDPP}.

%
%
The rest of the manuscript is organized as follows. First, in Section \ref{pdSec}, we introduce the concepts of $M$-bounded persistence diagrams and $M$-bounded persistent Betti numbers. Next, in Section \ref{resultSec}, we state the two main results of the paper, a CLT for the $M$-bounded persistence diagram and a functional CLT for the $M$-bounded persistent Betti numbers. In Section \ref{exSec}, we provide specific examples of point processes satisfying the conditions of the main results. Sections \ref{simSec} and \ref{dataSec} explore TDA-based tests for simulated and real datasets, respectively. Finally, Section \ref{discussSec} summarizes the findings and points to possible avenues of future research. The proofs of the main results are deferred to Sections \ref{proofSec} and \ref{verifSec} of the appendix.


\section{$M$-bounded persistent Betti numbers}
\label{pdSec}
%
%
For a locally finite point set $\X \subset \R^2$, persistent Betti numbers provide refined measures for the amount of clusters and voids on varying length scales. More precisely, we let
 \begin{equation}\label{balls}
U_r(\X) = \bigcup_{x\in \X} B_r(x).
\end{equation}	
denote the union of closed disks of radius $r \ge 0$ centered at points in $\X$. A 0-dimensional topological feature is a connected component of this union, corresponding to a cluster of points in $\X$, while a 1-dimensional feature can be thought of as a bounded connected component of the background space, often identified with its boundary loop, and describes a vacant area in the plane. As the disks grow, new features arise and vanish; we say that they are born and die again.
The persistent Betti numbers quantify this evolution of clusters and loops. Henceforth, we consider the persistence diagram only until a fixed deterministic radius $\tf \ge 0$.

%
%
As $r$ approaches the critical radius for continuum percolation, long-range phenomena emerge \cite{cPerc}. Thus, determining whether two points are connected could require exploring large regions in space. While useful quantitative bounds on cluster sizes are known for Poisson point processes \cite{teix}, for more general classes of point processes the picture remains opaque and research is currently at a very early stage \cite{gibbsPerc2, bartekPerc1}. Recently, a central limit theorem for persistent Betti numbers has been established in the Poisson setting \cite{krebs, shirai}, but for general point processes the long-range interactions pose a formidable obstacle towards proving a fully-fledged functional CLT.

From a more practical point of view, these long-range dependencies are of less concern. Although large features can carry interesting information, we expect that spatially bounded topological features already provide a versatile tool for the statistical analysis of both simulated point patterns and real datasets, even when focusing only on features of a bounded size. For that purpose, we concentrate on features whose spatial diameter does not exceed a large deterministic threshold $M$.

To define these \emph{$M$-bounded features}, we introduce the Gilbert graph $G_r(\X)$ on the vertex set $\X$. The Gilbert graph $G_r(\X)$ has for vertices the points in $\X$ and two points are connected by an edge if the distance between them is at most $2r$ or, equivalently, if the two disks of radius $r$ centered at the points intersect.

\subsection{$M$-bounded clusters}
The 0-dimensional $M$-bounded features alive at time $r > 0$ are the connected components of $G_r(\X)$ with diameter at most $M$. Starting at $r = 0$, all points belong to separate connected components that merge into larger clusters when $r$ increases. We thus say that all components are \emph{born at time 0}.

%
%
To define the death time of a component, let $\CO_r(x)$ denote the connected component of $x \in \X$ in $G_r(\X)$. The components of $x, y \in \X$ meet at time
$$R(x, y) = \inf\{r > 0:\, \CO_r(\x) = \CO_r(y)\}.$$
Then, the \emph{death time} of $\x \in \X$ is the smallest $R(x,y)$ such that the spatial diameter of $\CO_r(\x)$ exceeds $M$ or such that $P_x$ is lexicographically larger than $P_y$, where $P_x, P_y$ are the points of $\CO_r(x) \cap \X$ and $\CO_r(y) \cap \X$ whose associated disks meet at time $R(x, y)$. This ordering determines which component dies when two of them meet.

\subsection{$M$-bounded loops}
\label{loopSec}

Next, we introduce 1-dimensional features. At time $r > 0$, these correspond to holes, i.e., bounded connected components in the vacant phase $\V_r(\X) = \R^2 \setminus U_r(\X)$. 
In contrast to the clusters, there are no holes at time $0$, so that both birth and death times must be specified. Moreover, it needs to be defined how holes are related for different radii $r$. 
 
The \emph{death time} of a hole $H_s$ in $V_s(\X)$ is the first time $r > s$ when the hole is completely covered by disks, i.e., $H_s\subseteq U_r(\X)$. We identify a hole $H$ with the point $p(H)$ that is covered last. Thus, holes $H_s$ in $V_s(\X)$ and $H_r$ in $V_r(\X)$, are identified if $p(H_r)= p(H_s)$. 
 
New holes in $V_r(\X)$ can only appear when two balls merge, which corresponds to including a new edge in $G_r(\X)$.  When a new hole is formed, it can happen in two ways: either a finite component is separated from the infinite component, or an existing hole is split in two. In both cases, we define the size of the newly created hole(s) as follows: Let $x_1, \dots, x_k \in \X $ be the points in $\X$ such that the disks of radius $r$ around the points intersect the boundary of the hole $H$ in $V_r(\X)$. Then, \emph{the size of $H$} is the diameter of the set $\{x_1,\dots,x_k\} $. The size remains unchanged until the next time the hole is split into smaller pieces. Then the size is recomputed for both new holes. This definition ensures that the size decreases when the balls grow and only changes when a new edge is added to $G_r(\X)$.
 
The \emph{birth time} of a hole $H$ is the minimal $s$ such that there is a hole $H_s$ in $V_s(\X)$ with $p(H) = p(H_s)$ and size less than $M$. By an $M$-bounded loop, we mean a loop with size lower than $M$.

\subsection{The persistence diagram}
%
%
We now adapt the definition of the persistence diagram in \cite{shirai} to only include $M$-bounded features. That is, we define the \emph{$q$th $M$-bounded persistence diagram}, $q \in \{0, 1\}$, as the empirical measure
\begin{align}
\label{pdEq}
\PDD(\X) = \sum_{i \in I^{M,q}(\X)} \delta_{(B_i^M, D_i^M)},
\end{align}
where $I^{M, q}(\X)$ is an index set over all $M$-bounded $q$-dimensional features that die before time $\tf$ and $B_i^M, D_i^M$ are the birth and death times of the $i$th feature. Then, the \emph{$q$th $M$-bounded persistent Betti numbers} 
$$\b^{M, q}_{b, d}(\X) = \PDD(\X)([0, b] \times [d, \tf])$$
are the number of $M$-bounded features born before time $b \ge 0$ and dead after time $d \le \tf$.  When $q = 0$, all features are born at time 0, so that only death times are relevant. Hence, we write $\b^{M, 0}_d$ instead of the more verbose $\b^{M, 0}_{b, d}$.

\section{Main results}
\label{resultSec}
%
%
Henceforth, $\PP$ denotes a simple stationary point process in $\R^2$. We think of $\PP$ as a random variable taking values in the space of locally finite subsets $\NN$ of $\R^2$ endowed with the smallest $\sigma$-algebra $\mathfrak N$ such that the number of points in any given Borel set becomes measurable. Throughout the manuscript, we assume that the factorial moment measures exist and are absolutely continuous. In particular, writing  $\xx = (x_1, \dots, x_p) \in \R^{2p}$, the \emph{$p$th factorial moment density} $\rho^{(p)}$ is determined via the identity
\begin{align}
	\label{factMomEq}
	\E\Big[\prod_{i \le p} \PP(A_i)\Big] = \int_{A_1 \times \cdots \times A_p} \rho^{(p)}(\xx)\d \xx
\end{align}
for any pairwise disjoint bounded Borel sets $A_1, \dots, A_p \subset \R^2$, where $\PP(A_i)$ denotes the number of points of $\PP$ in $A_i$. Moreover, as we rely on the framework of \cite{yogeshCLT}, we also require that $\PP$ exhibits \emph{exponential decay of correlations}. Loosely speaking this expresses an approximate factorization of the factorial moment densities and is made precise in Section \ref{exSec} below. Many of the most prominent examples of point processes appearing in spatial statistics exhibit exponential decay of correlations \cite[Section 2.2]{yogeshCLT}.

%
%
Our first main result is a CLT for the persistence diagram built on the restriction $\PP_n = \PP \cap \Wn$ of the point process $\PP$ to a large observation window $\Wn = [-\sqrt n/2, \sqrt n/2]^2$. With a slight abuse of notation, we write $\PP \cup \xx = \PP \cup \{x_1, \dots, x_p\}$. 
To prove the CLT, we impose an additional condition concerning moments under the \emph{reduced $p$-point Palm distribution} $\P_{\xx}^!$. We recall that this distribution is determined via 
\begin{align}
	\label{redPalmEq}
	\E\Big[\sum_{{(X_1, \dots, X_p) \in \PP_{\ne}^p}} f(X_1, \dots, X_p; \PP) \Big] = \int_{\R^{2p}}\E_{\xx}^![f(\xx; \PP \cup \xx)]  \rho^{(p)}(\xx) \d \xx,
\end{align}
for any bounded measurable $f: \R^{2p} \times \NN \to \R $, where $\PP_{\ne}^p$ denotes $p$-tuples of pairwise distinct points in $\PP$. Then, we impose the following moment condition.
\begin{itemize}
	\item[{\bf (M)}] For every $p\ge1$ 
		$$\sup_{\substack{l \le p \\ \xx \in \R^{2l}}}\E_{\xx}^![\PP(\Wo)^p] < \infty.$$
\end{itemize}

%
%

%

To state the CLT for the persistence diagram precisely, we let 
$$\lan f, \PDDn \ran = \int_{\tfs}f(b, d) \PDD(\PP_n) (\d b, \d d) = \sum_{i \in I^{M, q}(\PP_n)} f\big(B_i^M, D_i^M \big)$$
denote the integral of a bounded measurable function $f:\,\tfs \to \R$ with respect to the measure $\PDDn$.

\begin{theorem}[CLT for persistence diagrams]
	\label{pdThm}
	Let $M > 0$, $q \in \{0, 1\}$ and $f:\,\tfs \to \R$ be a bounded measurable function. Assume that $\PP$ exhibits exponential decay of correlations and satisfies condition {\bf (M)}. Furthermore, assume that $\liminf_{n \to \infty}\Var(\lan f, \PDDn\ran) n^{-\nu} = \infty$ for some $\nu > 0$. Then, 
		$$\frac{\lan f, \PDDn\ran - \E[\lan f, \PDDn\ran]}{\sqrt{\Var(\lan f, \PDDn\ran)}}$$
	converges in distribution to a standard normal random variable as $n \to \infty$.
\end{theorem}

%
%
In order to derive a functional CLT for the persistent Betti numbers, we add a further constraint on $\PP$, which is needed to establish a lower bound on the variance via a conditioning argument in the vein of \cite[Lemma 4.3]{gibbsCLT}. For this purpose, we consider a random measure $\La$, which is jointly stationary with $\PP$ and which we think of as capturing additional useful information on the dependence structure of $\PP$. For instance, if $\PP$ is a Cox point process, we choose $\La$ to be the random intensity measure. If $\PP$ is a Poisson cluster process, then $\La$ would describe the cluster centers. If the dependence structure is exceptionally simple, it is also possible to take $\La = 0$. The idea of using additional information is motivated from conditioning on the spatially refined information coming from the clan-of-ancestors construction in Gibbsian point processes \cite{gibbsCLT}.

The point process $\PP$ is \emph{conditionally $m$-dependent} if $\PP \cap A$ and $\PP \cap A'$ are conditionally independent given $\sigma(\La, \PP \cap {A''})$
for any bounded Borel sets $A, A', A'' \subset \R^2$ such that the distance between $A$ and $A'$ is larger than some $m > 0$. Here, $\sigma(\La, \PP \cap {A''})$ denote the $\sigma$-algebra generated by $\La$ and $\PP \cap {A''}$.

Finally, we impose an absolute continuity-type assumption on the Poisson point process in a fixed box with respect to $\PP$ when conditioned on $\La$ and the outside points. More precisely, we demand that there exists $\rac > 6M$ with the following property, where $\mc Q$ denotes a homogeneous Poisson point process in the window $W_{\rac^2}$.
\begin{enumerate}
	\item[{\bf (AC)}]  
		Let $E_1, E_2 \in \mathfrak N$ be such that $\min_{i \in \{1,2\}}\P(\mc Q \in E_i) > 0$. Then, 
		$$ \E\big[\min_{i \in \{1,2\}} \P\big( \PP_{\rac^2} \in E_i\,|\, \La , \PP \setminus W_{\rac^2}\big) \big] > 0.$$
\end{enumerate}
Although {\bf (AC)} appears technical, Section \ref{exSec} illustrates that it is tractable for many commonly used point processes.

Since the persistent Betti numbers exhibit jumps at the birth- and death times of features, we work in the Skorokhod topology \cite[Section 14]{billingsley}.  

%
%
\begin{theorem}[Functional CLT for persistent Betti numbers]
\label{mainThm}
Let $M > 0$ and $\PP$ be a conditionally $m$-dependent point process with exponential decay of correlations and satisfying conditions {\bf (M)} and {\bf (AC)}. Then, the following convergence statements hold true.
	\begin{enumerate}
		\item[\bf q=0.]
			The one-dimensional process 
			$$\big\{n^{-1/2}\big(\PDodn - \E[\PDodn]\big)\big\}_{d \le \tf}$$
			converges weakly in Skorokhod topology to a centered Gaussian process.
		\item[\bf q=1.]The two-dimensional process
			$$\big\{n^{-1/2}\big(\PDoodn - \E[\PDoodn]\big)\big\}_{b, d \le \tf}$$
			converges weakly in Skorokhod topology to a centered Gaussian process.

	\end{enumerate}
\end{theorem}

%
%
Additionally, \cite[Theorem 1.12]{yogeshCLT} implies convergence of the rescaled variances. While Theorem \ref{pdThm} is an adaptation of \cite{yogeshCLT}, Theorem \ref{mainThm} is much more delicate.
As an application of Theorem \ref{mainThm}, we obtain a functional CLT for the following two characteristics, which are modified variants of the \emph{accumulated persistence function} from \cite{biscioTDA}:
	$$\APFz_r( \PP_n) = \sum_{i \in I^{M, 0}(\PP_n) } D_i^M \, \one{\{D_i^M \leq r\}}$$ 
	and 
	$$\APFo_r( \PP_n) = \sum_{i \in  I^{M, 1}(\PP_n)} (D_i^M-B_i^M)\one{\{B_i^M \leq r\}}.$$
	\begin{corollary}[Functional CLT for the APF]
		\label{apfCor}
		Let $M > 0$ and $\PP$ be as in Theorem \ref{mainThm}. Then, both 
			$\big\{n^{-1/2}(\APFz_r(\PP_n) - \E[\APFz_r(\PP_n)])\big\}_{r \le \tf}$ and 
			$\big\{n^{-1/2}(\APFo_r(\PP_n) - \E[\APFo_r(\PP_n)])\big\}_{r \le \tf}$ converge to centered Gaussian processes.
	\end{corollary}

\section{Examples of point processes}
\label{exSec}
\newcommand{\nc}[1]{{\color[rgb]{0.8, 0.1, 0} #1}}
\newtheorem{Def}[theorem]{Definition}

In this section, we give examples of point processes which satisfy the assumptions of our main theorems. More precisely, we show that log-Gaussian Cox processes with compactly supported covariance functions and Mat\'ern cluster processes both satisfy the conditions of Theorems \ref{pdThm} and \ref{mainThm}. We also show that the Ginibre point process satisfies the conditions of Theorem \ref{pdThm}.

Conversely, we do not expect that hard-core point processes satisfy the functional central limit theorem in the generality of Theorem \ref{mainThm}. Indeed, hard-core conditions put a strict lower bound on the death time of clusters and the birth time of loops. However, we believe that suitable repulsive point processes, where the hard-core conditions only need to be imposed with a certain probability can be embedded in the framework of Theorem \ref{mainThm}.

We first recall the definition of exponential decay of correlations from \cite{yogeshCLT}. To this end, we define the \emph{separation distance} between $\xx = \{x_1, \dots, x_p\} \subset \R^2$ and $\xx' = \{x_{p + 1}, \dots, x_{p + q}\} \subset \R^2$ as in \cite[Formula (1.3)]{yogeshCLT} via
\begin{align}
  \label{sepDistEq}
  \dist(\xx, \xx') = \inf_{\substack{i \le p \\ j \le q}} |x_i - x_{p + j}|.
\end{align}

\begin{Def}
Let $\PP$ be a stationary point process in $\R^2$, such that the $k$-point correlation function $\rho^{(k)}$ exists for all $k \ge 1$.
Then, $\PP$ exhibits \emph{exponential decay of correlations} if there exist $a < 1$, $\phi: [0, \infty) \to [0, \infty)$ such that
\begin{enumerate}
	\item $\lim_{t \to \infty} t^n \phi(t) = 0$ for all $n \ge 1$,
  \item $\liminf_{t \to \infty} \log \phi(t) / t^b < 0$ for some $b > 0$, 
  \item
		$$|\rho^{(p + q)}(\xx \cup \xx') - \rho^{(p)}(\xx) \rho^{(q)}(\xx')| \le (p + q)^{a (p + q)} \phi(\dist(\xx, \xx'))$$
for any $\xx = \{x_1, \dots, x_p\}, \xx' = \{x_{p + 1}, \dots, x_{p + q}\} \subset \R^2$.
\end{enumerate}
\end{Def}

\subsection{Log-Gaussian Cox process}
Let $Y = \{Y(x)\}_{x\in \R^2}$ be a stationary Gaussian process with mean $\mu \in \R$ and covariance function $c(x, x') = c(x - x')$. Then, the random measure on $\R^2$ defined as $\mathbf \Lambda(B) = \int_B \exp(Y(x))\d x$, for any Borel subset $B \subset \R^2$ has moments of any order. Let $\PP$ be a Cox process with random intensity measure $\mathbf \Lambda$, referred to as a Log-Gaussian Cox process. 
By \cite[Equation (7)]{CMW2}, the factorial moment densities of $\PP$ are given by 
\begin{align*}
	\rho^{(j)}(u_1, \dots, u_j) = \exp\left( j\mu + \frac{jc(0)}2 
		 \right)
		\prod_{1\le i < i' \le j} \exp(c(u_i - u_{i'})).
\end{align*}
To apply Theorems \ref{pdThm} and \ref{mainThm}, we assume that $c$ is bounded and of compact support, which ensures that $\PP$ exhibits exponential decay of correlation. 

We show below that condition $\mathbf{(M)}$ is satisfied. Let $\xx = (x_1, \dots, x_l) \in \R^{2l}$. According to \cite[Theorem 1]{CMW}, the Log-Gaussian Cox process $\PP$ under the reduced Palm version is also a Log-Gaussian Cox process $\PP_{\xx}$ with underlying Gaussian process $Y_{\xx}(x) = Y(x)+\sum_{i \le l}c(x, x_i)$. According to \cite[Equation (5.4.5)]{daleyVJ}, 
\begin{equation*}
	\E^!_{\xx} [\PP(W_1)^p] = \E[\PP_{\xx}(W_1)^p] = \sum_{1 \le j \le p}\Delta_{j, l, p}\int_{W_1^j} \rho^{(j)}_{\xx}(u_1, \dots, u_j) \d u_1\cdots \d u_j
\end{equation*}
for suitable coefficients $\Delta_{j, l, p} \in \R$, where $\rho^{(j)}_{\xx}(u_1, \dots, u_j)$ denotes the $j$th factorial moment density with respect to $\PP_{\xx}$. Therefore, it is enough to prove that 
\[\sup_{\xx \in \R^{2l}}\int_{W_1^j} \rho^{(j)}_{\xx}(u_1, \dots, u_j) \d u_1\cdots \d u_j < \infty, \] 
for all $j, l \ge 1$. Now, Equation (8) in \cite{CMW2} gives that
	\begin{multline*}
 	\rho^{(j)}_{\xx}(u_1, \dots, u_j) 
		 = 
	 \exp\Big( j\mu + \frac{jc(0)}2 +
		\sum_{\substack{1 \le i \le j \\ 1 \le k \le l}} c(u_i, x_k) \Big)
		\prod_{1\le i < i' \le j} \exp(c(u_i - u_{i'})).
 	\end{multline*}
where the right-hand side is bounded as $\mu$ and $c$ are bounded independently of $\xx$. This verifies condition \textbf{(M)}.
 
Since conditionally on $\mathbf \Lambda$, the point process $\PP$ is a Poisson point process, the conditional $m$-dependence property holds with $\Lambda = \mathbf \Lambda$.

It remains to verify condition $\mathbf{(AC)}$. By \cite[Equation (6.2)]{rasmus}, conditionally on $\Lambda$, the distribution of the point process $\PP_{\rac^2}$ admits the density with respect to a homogeneous Poisson point process $\mc Q $ with intensity 1 in $W_{\rac^2}$ given by
\[f_{\mathbf \Lambda}(\phi) = \exp({|W_{\rac^2}|-\mathbf \Lambda(W_{\rac^2})})\prod_{x\in \phi}\exp(Y(x)), \]
where $\phi \in \mathfrak N$. In particular, $f_{\mathbf \Lambda}(\phi)$ is strictly positive for all $\phi$. Therefore, if $E_1, E_2$ are two events such that $\min_{i \in \{1, 2\}}\P(\mc Q \in E_i) > 0$, then $\P(\PP_{\rac^2}\in E_i\,|\,\Lambda = \mathbf \Lambda) > 0$. This verifies condition $\mathbf{(AC)}$. 

\subsection{Mat\'ern cluster process}
Let $\eta$ be a homogeneous Poisson point process in $\R^2$ with intensity $\g > 0$. Given a realization of $\eta$, we define a family of independent point processes $(\Phi_x)_{x\in \eta}$, where $\Phi_x$, $x\in \eta$, is a homogeneous Poisson point process with intensity 1 in the disk $B_R(x)$ of radius $R > 0$ centered at $x \in \R^2$. The point process $\PP = \bigcup_{x\in \eta}\Phi_x$ is referred to as a \emph{Mat\'ern cluster process}. Since $\PP$ is $2R$-dependent, it exhibits exponential decay of correlations.

Next, we verify condition \textbf{(M)}. For this purpose, we deduce from \cite[Section 5.3.2]{CMW} that a Mat\'ern cluster process is a Cox process whose random intensity measure $\mathbf \Lambda$ has as density the random field $(\lambda(x))_{x\in \R^2}$ given by 
\[\lambda(x) = \gamma \eta(B_R(x)).\] 
Now, let $\xx = (x_1, \dots, x_l) \in \R^{2l}$ and $p\ge 1$ be fixed. From \cite[Equations (19) and (20)]{CMW} we obtain that
\[\E^!_{\xx} [ \PP(W_1)^p ] = \frac1{\E\big[ \prod_{i \le l}\lambda(x_i) \big] }\cdot \E\big[ \PP(W_1)^p\prod_{i\le l}\lambda(x_i) \big].\]
Since $\eta(B_R(x))$ is increasing in $\eta$ for every $x \in \R^2$ in the sense of \cite{LP}, the Harris-FKG inequality \cite[Theorem 20.4]{LP} gives that 
\begin{align*}
\E\big[ \prod_{i \le l}\lambda(x_i) \big] \ge \prod_{i\le l} \E[\lambda(x_i)] = (\gamma\pi R^2)^l, 
\end{align*}
where we used that $\lambda(x_i) = \eta(B_R(x_i))$ is a Poisson random variable with parameter $\pi R^2$. 
In order to bound $\E\Big[ \PP(W_1)^p \prod_{i \le l}\lambda(x_i) \Big]$, we first apply the H\"older inequality and stationarity, to arrive at 
\begin{align*}
	\E\Big[\PP(W_1)^p\prod_{i \le l}\lambda(x_i) \Big] \le \E\big[ \PP(W_1)^{p (l + 1)}\big]^{1 / (l + 1)}\E[\la(o)^{l + 1}]^{l / (l + 1)}.
\end{align*}
First, $\E[\la(o)^{l + 1}] = \E[\eta(B_R(o))^{l + 1}]$ is finite since $\eta$ is a Poisson point process. For the remaining part, we note that $\PP(W_1)\le \sum_{y\in \eta\cap (W_1\oplus B_R(o))}\#\Phi_x$, where $W_1\oplus B_R(o) = \{x + y:\, x \in W_1,\, y \in B_R(o)\}$ denotes the Minkowski sum. Hence, 
\begin{align*}
	\E\big[ \PP(W_1)^{p (l + 1)}\big] &\le \E\Big[ \big( \sum_{x\in \eta\cap (W_1\oplus B_R(o))}\#\Phi_x\big)^{p (l + 1)}\Big]\\
& \le \E\Big[ \prod_{x\in \eta\cap (W_1\oplus B_R(o))}e^{p (l + 1)\#\Phi_x} \Big]\\
	& =  \exp\big(\g |W_1\oplus B_R(o)| (\E[e^{p (l + 1) \#\Phi_0}] - 1) \big).
\end{align*}
where $\Phi_0$ is a homogeneous Poisson point process of intensity 1 in the disk $B_R(o)$ \cite[Theorem 3.9]{LP}. Again, since $\#\Phi_0$ is a Poisson random variable with parameter $\pi R^2$, the latter expression is finite.
 Taking the supremum over all $\xx$ and all $l\le p$, this verifies condition \textbf{(M)}. The point process $\PP$ is also conditionally $m$-dependent, by taking $m = 2R$ and $\Lambda = \eta$. 

It remains to prove $\mathbf{(AC)}$. By \cite[Equation (6.2)]{rasmus}, conditional on $\Lambda = \eta$, the distribution of $\PP_{\rac^2}$ admits the density 
\[f_\eta(\phi) = \gamma \exp({|W_{\rac^2}|-\mathbf \Lambda(W_{\rac^2})})\prod_{x\in \phi}\eta(B_R(x))\] 
with respect to the distribution of a homogeneous Poisson point process.
Now, consider the event on the event
\[\mc E = \{W_{\rac^2}\subset \eta \oplus B_R(o)\},\] 
the density $f_\eta$ is positive. Therefore, if $E_1, E_2$ are such that $\min_{i \in \{1,2\}}\P(\mc Q \in E_i) > 0$, then almost surely 
\[\min_{i \in \{1, 2\}} \P(\PP_{\rac^2} \in E_i | \eta )\mathbf1_\mc E(\eta) > 0.\] 
Since $\mc E$ occurs with positive probability, this proves condition \textbf{(AC)}.

\subsection{Ginibre point process}
The Ginibre point process is a determinantal point process with kernel \[K(z_1, z_2) = \exp(z_1\overline{z_2})\exp\left(- \frac{|z_1|^2 + |z_2|^2}2 \right), \] with $z_1, z_2\in \mathbb C$. As mentioned in \cite[p. 19]{yogeshCLT}, this point process exhibits exponential decay of correlation. According to \cite[Theorem 2]{gold}, for $\xx = (x_1, \dots, x_l)\in \R^{2l}$ we have $\E^!_{\xx}\left[\PP(W_1)^p \right] \le \E{\left[\PP(W_1)^p\right]}$, where the right-hand side is finite by \cite[Lemma 4.2.6]{hough}. Hence, we obtain an upper bound for $\E^!_{\xx}\left[\PP(W_1)^p \right]$, which is independent of $\xx$, thereby verifying condition \textbf{(M)}.

\section{Simulation study}
\label{simSec}
We elucidate in a simulation study, how cluster- and loop-based test statistics derived from Theorem \ref{mainThm} can detect deviations from complete spatial randomness. The simulations are carried out on top of the \texttt R-packages \texttt{spatstat} and \texttt{TDA} \cite{fasy, spatstat}.

For the entire simulation study, the null model $\Poi(2)$ is a Poisson point process with intensity $2$ in a $10 \times 10$ observation window. Moreover, we fix $M = \sqrt 2 \cdot 10$ so large that it encompasses the entire sampling window and therefore suppress its appearance in the notation. Although the proof of Theorem \ref{mainThm} relies on the $M$-boundedness, the simulation study illustrates that it is not critical to impose this condition when testing hypotheses on common point patterns.

\subsection{Deviation tests}
\label{devTestSec}
As a first step, we derive scalar cluster- and loop-based test statistics.

%
%
\subsubsection{Definition of test statistics}

%
%
As a test statistic based on clusters, we use the integral over the number of cluster deaths in a time interval $[0, \rc]$ with $\rc \le \tf$, i.e.,  
\begin{align}
	\label{t0Def}
	\Tc = \int_0^{\rc} \ms{PD}^0(\PPn)([0, d]) {\rm d}d.
\end{align}
After subtracting the mean, this test statistic becomes reminiscent of the classical Cram\'er-von-Mises statistic except that we do not consider squared deviations. Although squaring would make it easier to detect two-sided deviations, it would also require knowledge of quantiles of the square integral of a centered Gaussian process. Albeit possible, this incurs substantial computational expenses. Our simpler alternative has the appeal that as an integral of a Gaussian process, $\Tc$ is asymptotically normal and therefore characterized by its mean and variance.

%
%
As a test statistic based on loops, we use the accumulated persistence function, which aggregates the life times of all loops with birth times in a time interval $[0, \rl]$ with $\rl \le \tf$, i.e.,
\begin{align}
	\label{t1Def}
	\Tl = \ms{APF}^1_{\rl}(\PPn) = \int_{[0, \rl] \times [0, \tf]} (d - b) \ms{PD}^1(\PPn)({\rm d}b, {\rm d}d).
\end{align}
 By Corollary \ref{apfCor},  after centering and rescaling, the statistic $\Tl$ converges in the large-volume limit to a normal random variable.

 The statistics $\Tc$ and $\Tl$ are specific possibilities to define scalar characteristics from the persistence diagram. Depending on the application context other choices, such as $\ms{APF}^0$ instead of $\Tc$ could be useful. However, in the simulation study below we found the weighting by life times of clusters to be detrimental.

\subsubsection{Exploratory Analysis}
\label{altSec}
\begin{figure}[!htpb]
	\includegraphics[trim={0 10cm 0 10cm},clip,width=0.9\textwidth]{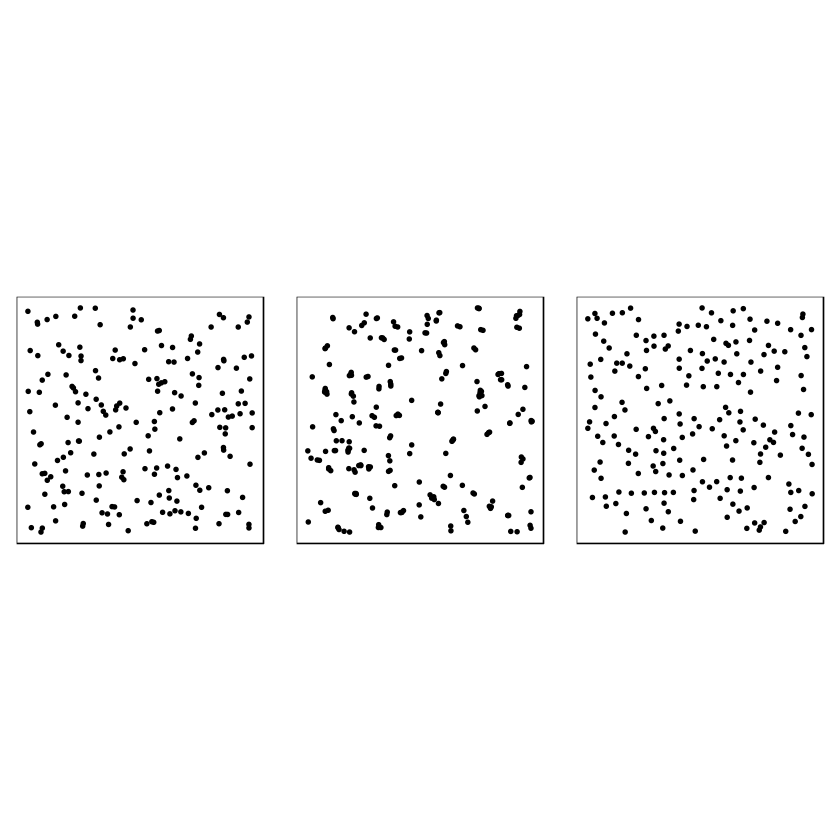}
	\caption{Samples from the $\Poi(2)$ null model (left), the $\mcp(2, 0.1, 1)$ process (center) and the $\str(4.5, 0.1, 0.35)$ process (right).}
	\label{ppFig}
\end{figure}

As alternatives to the Poisson null hypothesis, we consider the attractive Mat\'ern cluster and the repulsive Strauss processes. More precisely, the Mat\'ern cluster process $\mcp(2, 0.1, 1)$ features a Poisson parent process with intensity 2 and generates a $\Poi(1)$ number of offspring uniformly in a disk of radius $0.1$ around each parent. The Strauss process $\str(4.5, 0.1, 0.35)$ has interaction parameter $0.1$ and interaction radius $0.35$. The intensity parameter $4.5$ was tuned so as to match approximately the intensity of the null model. Figure \ref{ppFig} shows realizations of the null model and the alternatives.

\begin{figure}[!htpb]
	\includegraphics[trim={0 3.2cm 0 3.1cm}, clip, width=.85\textwidth]{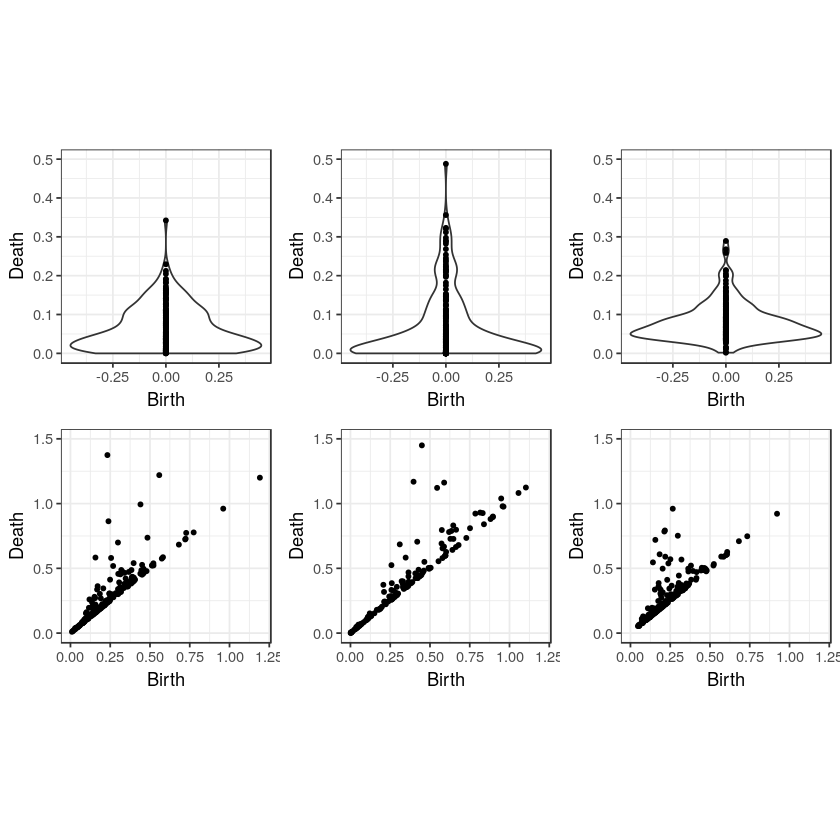}
	\caption{Persistence diagrams for cluster-based features with density plots (top) and loop-based features (bottom) for the $\Poi(2)$ null model (left), the $\mcp(2, 0.1, 1)$ process (center) and the $\str(4.5, 0.1, 0.35)$ process (right).}
	\label{pdSimFig}
\end{figure}

In a first step, in Figure \ref{pdSimFig}, we plot the persistence diagrams of samples from the null model and of the alternatives.

From the cluster-based diagrams, it becomes apparent that in comparison to the null model, in the Mat\'ern cluster process, features can die also at rather late times, whereas this happens very rarely in the Strauss process. When analyzing loops, we see that loops with long life times can appear earlier in the null model than in the Mat\'ern cluster process. Conversely, while some loops with substantial life time emerge at later times  in the null model, there are very few such cases in the Strauss model.

\subsubsection{Mean and variance under the null model}
Now, we determine the mean and variance of $\Tc$ and $\Tl$ under the null model with $\tf = 1.5$. For this purpose, we compute the number of cluster deaths and accumulated loop life times for 10,000 independent draws of the null model.

\begin{figure}[!htpb]
	\includegraphics[width=0.55\textwidth]{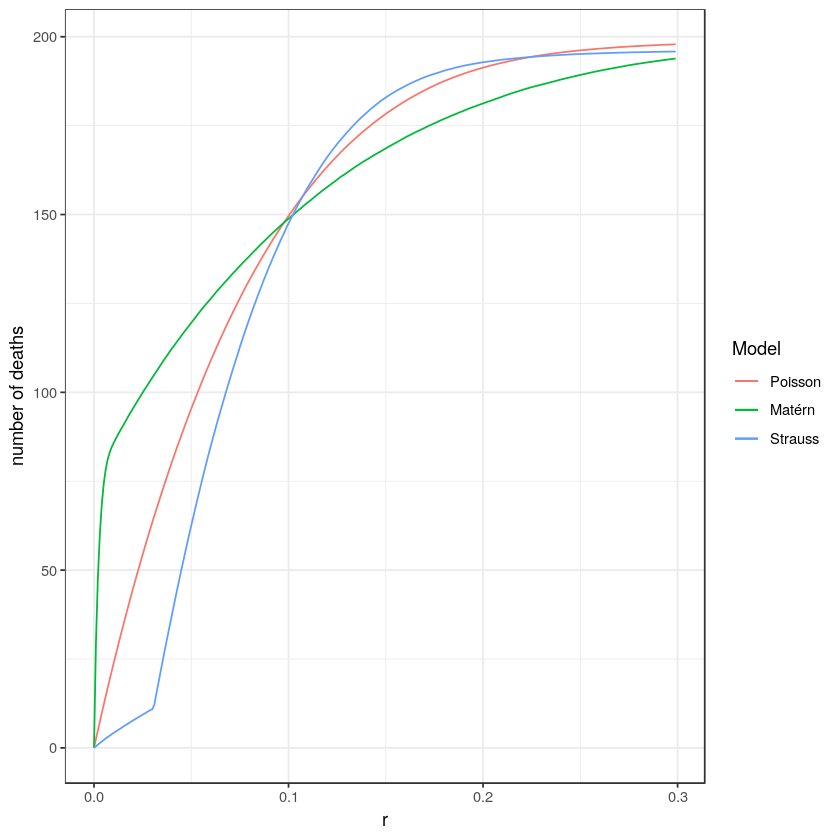}~\includegraphics[width=0.55\textwidth]{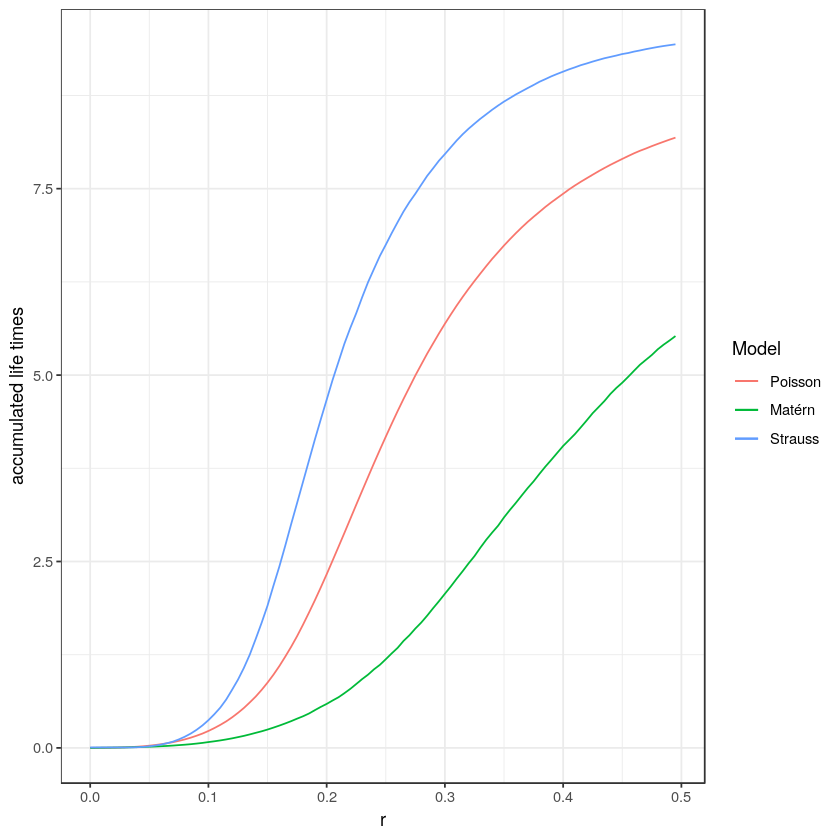}
	\caption{Mean number of cluster deaths (left) and accumulated loop life times (right) for the null model (red) and the alternatives (green and blue) based on 10,000 realizations.}
	\label{meanCurveFig}
\end{figure}

Comparing the mean curves for the number of cluster deaths in the null model with those of the alternatives matches up nicely with the intuition about attraction and repulsion. For late times, they all approach a common value, namely the expected number of points in the observation window. However, Figure \ref{meanCurveFig} shows that for the Mat\'ern model, the slope is far steeper for early times, caused by merging of components of points within a cluster. In contrast, for the Strauss process the increase is at first much less pronounced than in the Poisson model, thereby reflecting the repulsive nature of the Gibbs potential.

For the loops, a radically different picture emerges. Here, the curve for the Strauss process lies above the accumulated loop life times of the null model. The Strauss model spawns substantially more loops than the Poisson model, although most of them live for a shorter period. Still, taken together these competing effects lead to a net increase of the accumulated loop life times in the Strauss model.
\vspace{-.3cm}

\subsubsection{Type I and II errors}
By Theorem \ref{mainThm}, the statistics $\Tc$ and $\Tl$ are asymptotically normal, so that knowing the mean and variance allows us to construct a deviation test whose nominal confidence level is asymptotically exact. For the loops, we can choose the entire relevant time range, so that  $\rl = 0.5$. For the cluster features, this choice would be unreasonable, as for late times, we simply obtain the number of points in the observation window, which is not discriminative. Hence, we set $\rc = 0.1$. We stress that in situations with no a priori knowledge of a good choice of $\rc$, the test power can degrade substantially.

To analyze the type I and II errors, we draw 1,000 realizations from the null model and from the alternatives, respectively. Table \ref{powTab} shows the rejection rates of this test setup. Under the null model the rejection rates are close to the nominal 5\%-level, thereby illustrating that already for moderately large point patterns the approximation by the Gaussian limit is accurate. Using the mean and variance from the null model, we now compute the test powers for the alternatives. Already $\Tc$ leads to a test power of approximately $60\%$ for both alternatives. When considering $\Tl$, we obtain a type I error rate of 4.8\%, so that the confidence level is kept. Moreover, the power analysis reveals that in the present simulation set-up, $\Tl$ is better in detecting deviations from the null hypothesis than $\Tc$.

	 \begin{table}[!htpb]
 \begin{tabular}{lllr}
	 & $\Poi$ & $\mcp$ & $\str$  \\
\hline
	 $\Tc$	       & 5.1\%    & 59.3\% & 60.7\%       \\
	 $\Tl$	        & 4.8\% & 94.7\%  &   71.4\%  \\ 
\end{tabular}
	 \caption{Rejection rates for the test statistics $\Tc$ and $\Tl$ under the null model and the alternatives.}
	 \label{powTab}
	 \end{table}

\subsection{Envelope Tests}

Leveraging Theorem \ref{mainThm} shows that the deviation statistics $\Tc$ and $\Tl$ are asymptotically normal. Using a simulation-based estimate for the asymptotic mean and variance under the null model allowed us to construct a deviation test whose confidence level is asymptotically precise. A caveat of the above analysis is that the magnitude of clustering and repulsion is strong and clearly visible in the samples.

Recently, global envelope tests have gained widespread popularity, because they are both powerful and provide graphical insights as to why a null hypothesis is rejected \cite{envelop}. The global envelope tests are fundamentally Monte Carlo-based tests and therefore do not relate directly to the large-volume CLT. However, they also rely on a functional summary statistics as input. Most of the applications in spatial statistics use a distance-based second-order functional such as Ripley's $L$-function. In this section, we compare such classical choices with cluster- and loop-based statistics.

\subsubsection{Alternatives}
Since envelope tests excel at detecting subtle changes from the null model, we consider now a new parameter set-up to compare the $L$-function with the cluster- and loop-based statistics. Here, both the Mat\'ern cluster as well as the Strauss process are substantially more similar to the Poisson point process. Hence, for the alternatives, we use again Mat\'ern cluster and Strauss processes, but choose different parameters.

\begin{figure}[!htpb]
	\includegraphics[trim={0 9cm 0 0cm},clip,width=0.9\textwidth]{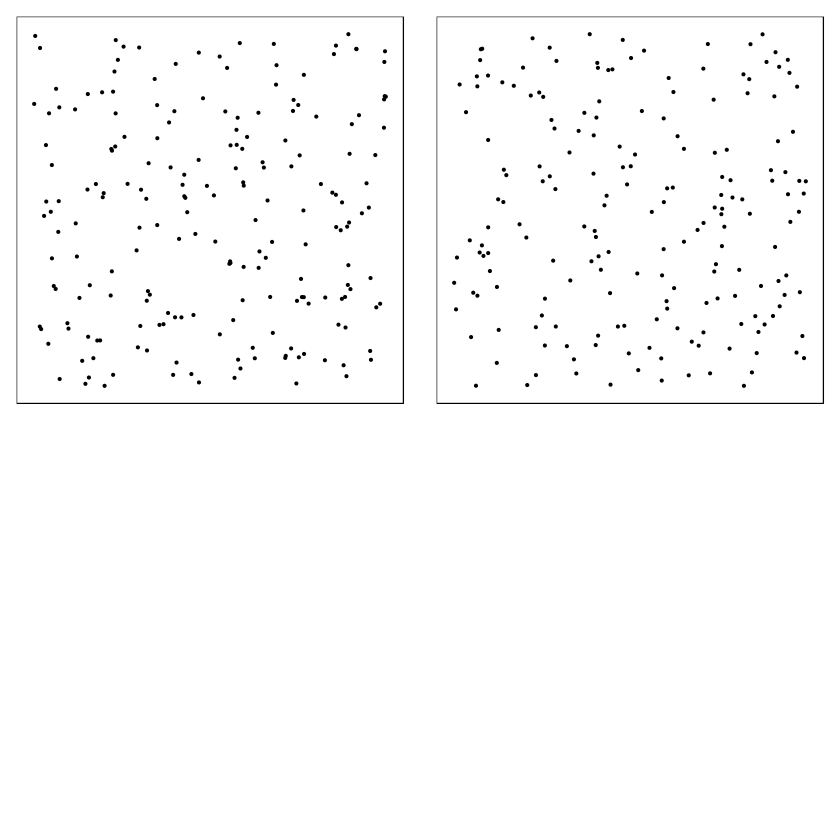}
	\caption{Samples from $\mcp(20, 0.1, 0.1)$ (left) and $\str(2.1, 0.1, 0.1)$ (right).}
	\label{ppsEnvFig}
\end{figure}

We found that the cluster- and loop-based statistics were particularly powerful in situations involving small interaction radii. Hence, as alternatives we choose the $\mcp(20, 0.1, 0.1)$ process and the $\str(2.1, 0.1, 0.1)$ process, see Figure \ref{ppsEnvFig}. The interaction parameter of the Strauss process was again tuned to match approximately the intensity of the null model. 

\newpage
\subsubsection{Power analysis}
To analyze the power of the envelope test, we generate $s = 4,999$ realizations of the null model and 1,000 realizations of the alternatives. Then, we perform the global envelope test from \cite{envelop} with three functional summary statistics. The first is Ripley's $L$-function \cite[Definition 4.6]{rasmus}. Second, we consider the number of cluster deaths as illustrated in Figure \ref{meanCurveFig}. Third, for the loops, we use a two-dimensional functional statistics derived from the persistent Betti numbers $\{\beta^{*, 1}_{b, l}\}_{b, l}$ associated with life times $l$ rather than death times $d$ in order to expand the support of the statistic to the entire first quadrant. More precisely, $\beta^{*, 1}_{b, l}$ counts the number of loops born before time $b$ and with life time at least $l$. 

The rejection rates from Table \ref{powEnvTab} illustrate that for the alternatives described above, the cluster-based test gives a similar test power as the $L$-function-based test for the Mat\'ern cluster process and a substantially higher test power for the Strauss process. Moreover, the loop-based test works even better in the Strauss case, but performs substantially worse for the Mat\'ern alternative.

 \begin{center}
	 \begin{table}[!htpb]
 \begin{tabular}{lllr}
	 & $\mcp$ & $\str$  \\
\hline
	 Ripley's $L$	           & 42.6\% & 20.5\%       \\
	 $\ms{Cluster}$	         & 41.5\%  &   26.3\%  \\ 
	 $\ms{Loop}$	         & 27.0\%  &   32.2\%  \\ [.5ex]
\end{tabular}
	 \caption{Rejection rates for envelope tests based on the $L$ function and cluster- and loop-based functional statistics.}
	 \label{powEnvTab}
	 \end{table}
 \end{center}

\section{Analysis of the minicolumn dataset}
\label{dataSec}
In this section, we explore to what extent the deviation tests from Section \ref{simSec} provide insights when dealing with real data. For this purpose, we analyze the minicolumn dataset provided by scientists at the Centre for Stochastic Geometry and Advanced Bioimaging. 

As it should serve only to illustrate the application of Theorem \ref{mainThm}, the present analysis is very limited in scope, and we refer to \cite{mc3} for a far more encompassing study. For instance, that work considers two datasets and investigates 3D data together with marks for the directions attached to the neurons. 

\subsection{Exploratory analysis}
\begin{figure}[!htpb]
	\includegraphics[width=0.52\textwidth]{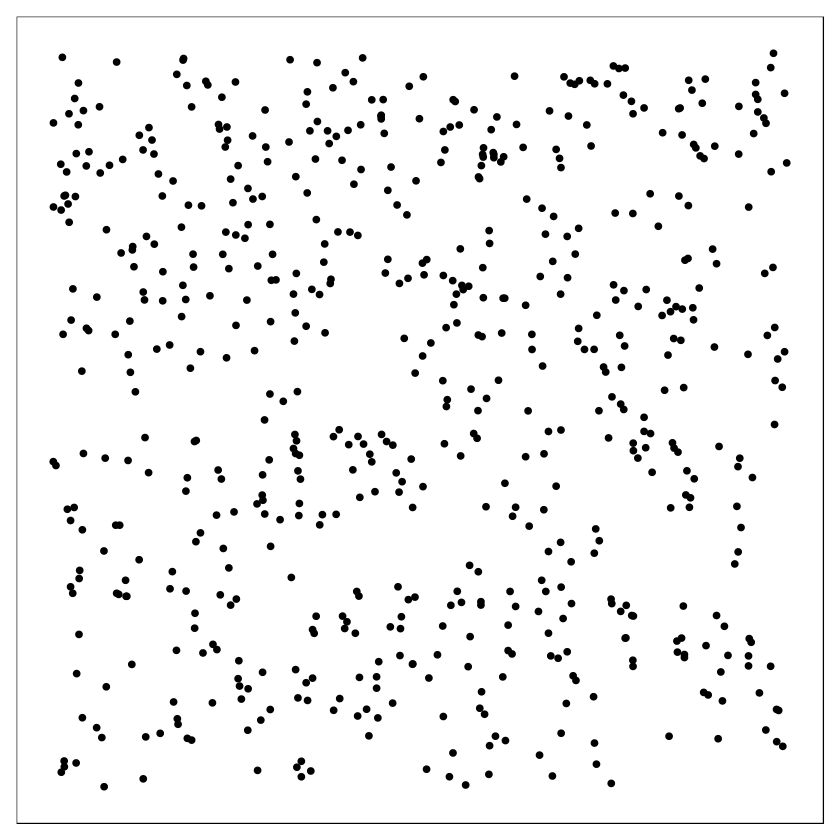}
	\caption{Projected minicolumn point pattern}
	\label{miniPatFig}
\end{figure}

The minicolumn dataset consists of 634 points emerging as two-dimensional projections of a three-dimensional point pattern of neurons. As neurons are believed to arrange in vertical columns, the projections are expected to exhibit clustering, see \cite{mc1, mc2}. The projections are taken along $z$-axis, since neuroscientists expect an arrangement in vertical columns. A visual inspection of the point pattern in Figure \ref{miniPatFig} supports this hypothesis. 

\begin{figure}[!htpb]
	\includegraphics[trim={0 4cm 0 4cm},clip, width=0.8\textwidth]{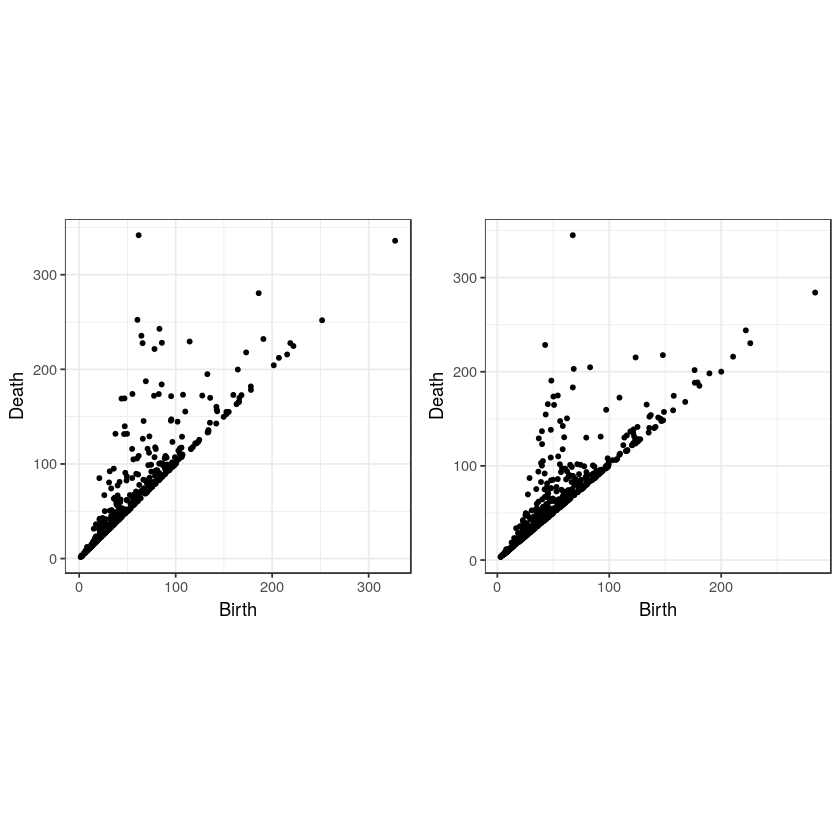}
	\caption{Persistence diagram for the minicolumn data (left) and a homogeneous Poisson point process with the same intensity (right)}
	\label{miniPDFig}
\end{figure}

As a first step, we explore whether the purported clustering already manifests in the persistence diagram. Comparing the loop-based persistence diagram of the minicolumn data with the persistence diagram of a homogeneous Poisson point process in Figure \ref{miniPDFig} shows that loops with substantial life times tend to be born later in the minicolumn model. This suggests clustering since loops formed by points within a cluster typically disappear rapidly.

\begin{figure}[!htpb]
	\includegraphics[width=0.48\textwidth]{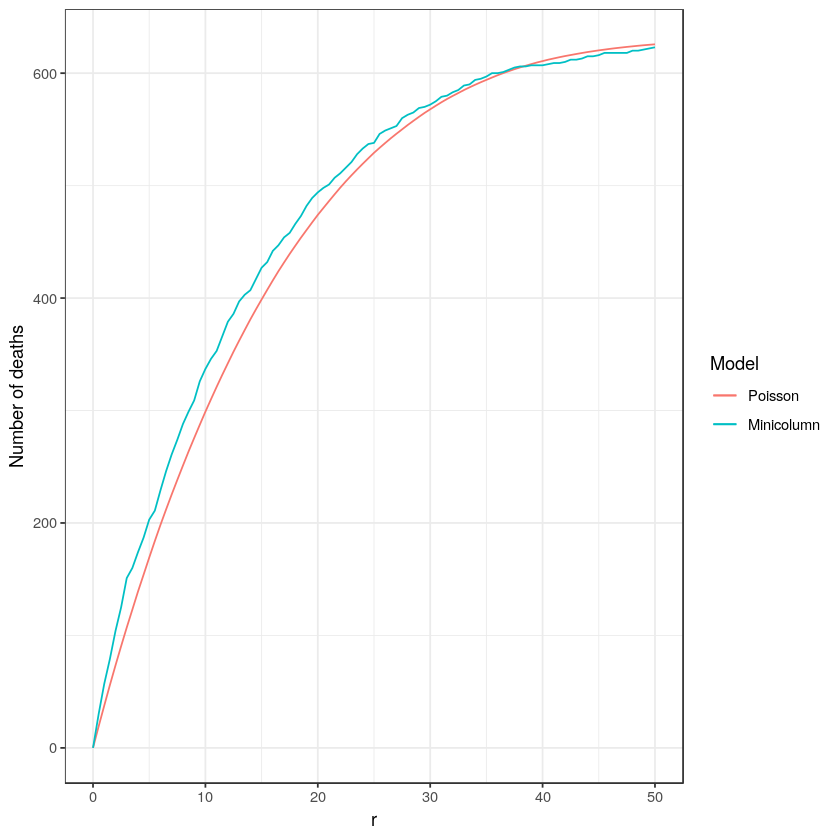}\quad
	\includegraphics[width=0.48\textwidth]{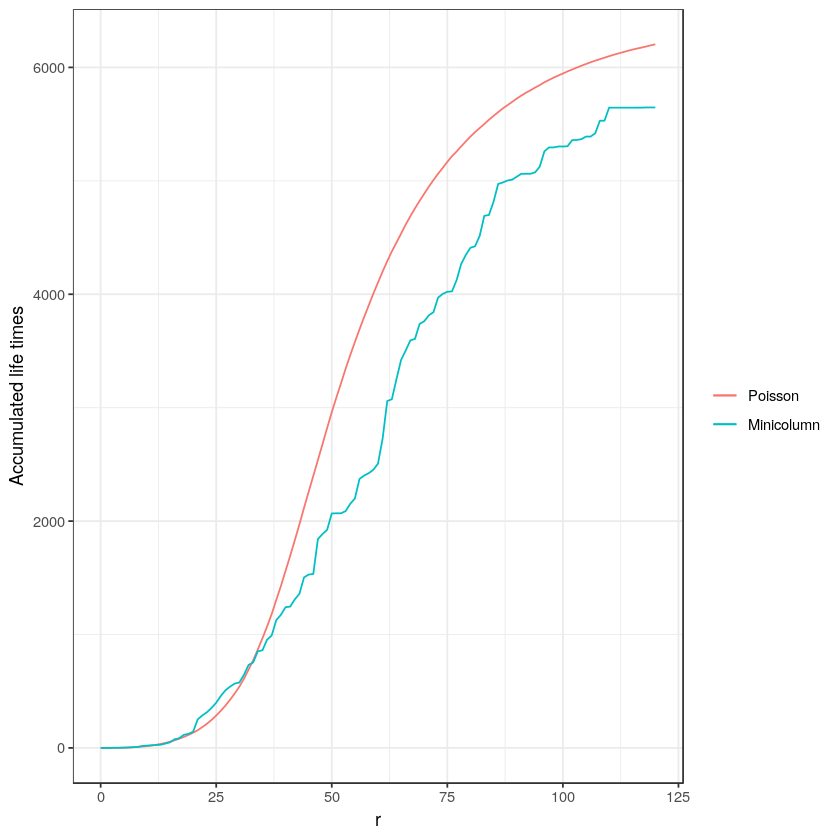}
	\caption{Number of cluster deaths (left) and accumulated loop life times (right) for the Poisson null model (red) and the minicolumn dataset (green).}
	\label{minicurvesFig}
\end{figure}

Now, we explore whether the impressions from the persistence diagrams are reflected in the summary statistics from Section \ref{simSec}. When comparing in Figure \ref{minicurvesFig} (left) the number of cluster death at different points in time, we note that until time 35, the curve for the observed data runs a bit above the curve for the null model. This provides already a first indication towards clustering. Next, we proceed to the loop-based features. As shown in Figure \ref{minicurvesFig} (right), the curve for the observed pattern runs substantially below the one of the null model. This reflects a property that we have seen already in the persistence diagram: clusters with substantial life time tend to be born earlier in the null model, thereby leading to a steeper increase of the accumulated life times.

\subsection{Test for complete spatial randomness}
Under the impression of the previous visualizations, we now test the minicolumn pattern against the null model. As in Section \ref{simSec}, we deduce from Theorem \ref{mainThm} that the statistics are asymptotically normal under the null model, so that we only need to determine means and variances. 

A subtle issue concerns the choice of the integration interval. The simplest option would be to take the whole intervals shown in Figure \ref{minicurvesFig}. For instance, for the loop-based features, this means $\rl = \tf = 120$. However, for the cluster-based features the choice of the interval is less clear, since taking the whole interval is not discriminatory. The experiences from the simulation study indicate that the test is most powerful for early death times. Therefore, we choose $\rc = 10$.

With these choices, both the cluster-based and the loop-based test reject the null-hypothesis at the 5\% level, since the corresponding $p$-values are $1.7\%$ and $1.2\%$. However, the tests are sensitive to the choice of the integration bound. This is especially true for the cluster-based test, where going from $\rc = 10$ to $\rc = 15$ results in a $p$-value of $5.2\%$, so that the null hypothesis is no longer rejected. On the other hand, the loop-based test still yields a $p$-value of $1.9\%$ $\rl = 75$. However, when reducing even further to $\rl = 50$, then the $p$-value increases sharply to $34.1\%$, so that the null-hypothesis is no longer rejected.

\section{Discussion}
\label{discussSec}
%
%
In this paper, we elucidated how to apply tools from TDA to derive goodness-of-fit tests for planar point patterns. For this purpose, we derived sufficient conditions for a large-domain functional CLT for the $M$-bounded persistent Betti numbers on point processes exhibiting exponential decay of correlations. Following the framework developed in \cite{yogeshCLT}, the main difficulty arose from a detailed analysis of geometric configurations when bounding higher-order cumulants.

%
%
A simulation study revealed that the asymptotic Gaussianity is already accurate for patterns consisting of a few hundred data points. Additionally, as functional summary statistics, the persistent Betti numbers can also be used in the context of global envelope tests. Here, our finding is that TDA-based statistics can provide helpful additional information for point patterns with small interaction radii.

%
%
Finally, we applied the TDA-based tests on a point pattern from a neuroscientific dataset. As conjectured from the application context, the functional summary statistics indicate a clustering of points and the tests reject the Poisson null-model. However, the analysis also reveals a sensitivity to the range of birth times considered in the statistics.

%
%
In future work, we plan to extend the present analysis to dimensions larger than 2. On a technical level, the definition of higher-dimensional features requires a deeper understanding of persistent homology groups. Additionally, when thinking of broader application scenarios, a further step is to extend the testing framework from mere point patterns to random closed sets involving a richer geometric structure.

\section*{Acknowledgments}
We thank J.~M\o ller and A.~D.~Christoffersen for valuable discussions on the minicolumn dataset and for providing references. We also thank J.~Yukich for inspiring discussions and helpful remarks. Finally, we thank the Centre for Stochastic Geometry and Advanced Bioimaging for collecting and sharing the data. CB, CH and AS are supported by The Danish Council for Independent Research | Natural Sciences, grant DFF – 7014-00074 \emph{Statistics for point processes in space and beyond}, and by the \emph{Centre for Stochastic Geometry and Advanced Bioimaging}, funded by grant 8721 from the Villum Foundation. NS was partially supported by the French ANR grant ASPAG (ANR-17-CE40- 0017).

\bibliography{./lit}

\begin{thebibliography}{10}

\bibitem{teix}
D.~Ahlberg, V.~Tassion, and A.~Teixeira.
\newblock Sharpness of the phase transition for continuum percolation in
  $\mathbb {R}^2$.
\newblock {\em Probab. Theory Related Fields}, 172(1):525--581, 2018.

\bibitem{spatstat}
A.~Baddeley and R.~Turner.
\newblock spatstat: An {R} package for analyzing spatial point patterns.
\newblock {\em Journal of Statistical Software, Articles}, 12(6):1--42, 2005.

\bibitem{bary}
Y.~Baryshnikov and J.~E. Yukich.
\newblock Gaussian limits for random measures in geometric probability.
\newblock {\em Ann. Appl. Probab.}, 15(1A):213--253, 2005.

\bibitem{bickel}
P.~J. Bickel and M.~J. Wichura.
\newblock Convergence criteria for multiparameter stochastic processes and some
  applications.
\newblock {\em Ann. Math. Statist.}, 42:1656--1670, 1971.

\bibitem{billingsley}
P.~Billingsley.
\newblock {\em Convergence of Probability Measures}.
\newblock J. Wiley \& Sons, New York, second edition, 1999.

\bibitem{biscioTDA}
C.~A.~N. Biscio and J.~M{\o}ller.
\newblock The accumulated persistence function, a new useful functional summary
  statistic for topological data analysis, with a view to brain artery trees
  and spatial point process applications.
\newblock {\em J. Comput. Graph. Statist.}, 2019 (to appear).

\bibitem{bartekPerc1}
B.~B{\l}aszczyszyn and D.~Yogeshwaran.
\newblock Clustering and percolation of point processes.
\newblock {\em Electron. J. Probab.}, 18:1--20, 2013.

\bibitem{yogeshCLT}
B.~B{\l}aszczyszyn, D.~Yogeshwaran, and J.~E. Yukich.
\newblock Limit theory for geometric statistics of point processes having fast
  decay of correlations.
\newblock {\em Ann. Probab.}, 47(2):835--895, 2019.

\bibitem{bubenik}
P.~Bubenik.
\newblock Statistical topological data analysis using persistence landscapes.
\newblock {\em J. Mach. Learn. Res.}, 16:77--102, 2015.

\bibitem{calka}
P.~Calka, T.~Schreiber, and J.~E. Yukich.
\newblock Brownian limits, local limits and variance asymptotics for convex
  hulls in the ball.
\newblock {\em Ann. Probab.}, 41(1):50--108, 2013.

\bibitem{carlsson}
G.~Carlsson.
\newblock Topology and data.
\newblock {\em Bull. Amer. Math. Soc.}, 46(2):255--308, 2009.

\bibitem{chazal}
F.~Chazal and M.~Bertrand.
\newblock High-dimensional topological data analysis.
\newblock In C.~D. Toth, J.~O'Rourke, and J.~E. Goodman, editors, {\em Handbook
  of Discrete and Computational Geometry}. CRC Press, Boca Raton, third
  edition, 2017.

\bibitem{silhouette}
F.~Chazal, B.~T. Fasy, F.~Lecci, A.~Rinaldo, and L.~Wasserman.
\newblock Stochastic convergence of persistence landscapes and silhouettes.
\newblock {\em J. Comput. Geom.}, 6(2):140--161, 2015.

\bibitem{mc3}
A.~D. Christoffersen, J.~M{\o}ller, and H.~S. Christensen.
\newblock Modelling columnarity of pyramidal cells in the human cerebral cortex
  based on directionally marked 3{D} point patterns. {I}n preparation.

\bibitem{CMW2}
J.-F. Coeurjolly, J.~M{\o}ller, and R.~Waagepetersen.
\newblock Palm distributions for log {G}aussian {C}ox processes.
\newblock {\em Scand. J. Stat.}, 44(1):192--203, 2017.

\bibitem{CMW}
J.-F. Coeurjolly, J.~M{\o}ller, and R.~Waagepetersen.
\newblock A tutorial on {P}alm distributions for spatial point processes.
\newblock {\em Int. Stat. Rev.}, 85(3):404--420, 2017.

\bibitem{daleyVJ}
D.~J. Daley and D.~Vere-Jones.
\newblock {\em An Introduction to the Theory of Point Processes}.
\newblock Springer-Verlag, New York, second edition, 2003.

\bibitem{edHar}
H.~Edelsbrunner and J.~Harer.
\newblock {\em Computational Topology}.
\newblock American Mathematical Society, Providence, RI, 2010.

\bibitem{raic3}
P.~Eichelsbacher, M.~Rai\v{c}, and T.~Schreiber.
\newblock Moderate deviations for stabilizing functionals in geometric
  probability.
\newblock {\em Ann. Inst. Henri Poincar\'e Probab. Stat.}, 51(1):89--128, 2015.

\bibitem{fasy}
B.~T. Fasy, J.~Kim, F.~Lecci, and C.~Maria.
\newblock Introduction to the {R} package {TDA}.
\newblock {\em arXiv preprint arXiv:1411.1830}, 2014.

\bibitem{gold}
A.~Goldman.
\newblock The {P}alm measure and the {V}oronoi tessellation for the {G}inibre
  process.
\newblock {\em Ann. Appl. Probab.}, 20(1):90--128, 2010.

\bibitem{kclt}
L.~Heinrich.
\newblock Gaussian limits of empirical multiparameter {$K$}-functions of
  homogeneous {P}oisson processes and tests for complete spatial randomness.
\newblock {\em Lith. Math. J.}, 55(1):72--90, 2015.

\bibitem{alphaDPP}
L.~Heinrich.
\newblock On the strong {B}rillinger-mixing property of
  {$\alpha$}-determinantal point processes and some applications.
\newblock {\em Appl. Math.}, 61(4):443--461, 2016.

\bibitem{heinrichShot}
L.~Heinrich and V.~Schmidt.
\newblock Normal convergence of multidimensional shot noise and rates of this
  convergence.
\newblock {\em Adv. in Appl. Probab.}, 17(4):709--730, 1985.

\bibitem{shirai}
Y.~Hiraoka, T.~Shirai, and K.~D. Trinh.
\newblock Limit theorems for persistence diagrams.
\newblock {\em Ann. Appl. Probab.}, 28(5):2740--2780, 2018.

\bibitem{hough}
J.~B. Hough, M.~Krishnapur, Y.~Peres, and B.~Vir\'{a}g.
\newblock {\em Zeros of {G}aussian Analytic Functions and Determinantal Point
  Processes}.
\newblock American Mathematical Society, Providence, 2009.

\bibitem{gibbsPerc2}
S.~Jansen.
\newblock Continuum percolation for {G}ibbsian point processes with attractive
  interactions.
\newblock {\em Electron. J. Probab.}, 21:No. 47, 22, 2016.

\bibitem{kallenberg}
O.~Kallenberg.
\newblock {\em Foundations of Modern Probability}.
\newblock Springer, New York, second edition, 2002.

\bibitem{krebs}
J.~T.~N. Krebs and W.~Polonik.
\newblock On the asymptotic normality of persistent {B}etti numbers.
\newblock {\em arXiv preprint arXiv:1903.03280}, 2019.

\bibitem{LP}
G.~Last and M.~Penrose.
\newblock {\em Lectures on the {P}oisson process}.
\newblock Cambridge University Press, Cambridge, 2018.

\bibitem{cPerc}
R.~Meester and R.~Roy.
\newblock {\em Continuum Percolation}.
\newblock Cambridge University Press, Cambridge, 1996.

\bibitem{mc1}
J.~M{\o}ller, F.~Safavimanesh, and J.~G. Rasmussen.
\newblock The cylindrical {$K$}-function and {P}oisson line cluster point
  processes.
\newblock {\em Biometrika}, 103(4):937--954, 2016.

\bibitem{rasmus}
J.~M{\o}ller and R.~P. Waagepetersen.
\newblock {\em Statistical Inference and Simulation for Spatial Point
  Processes}.
\newblock CRC, Boca Raton, 2004.

\bibitem{envelop}
M.~Myllym\"{a}ki, T.~Mrkvi\v{c}ka, P.~Grabarnik, H.~Seijo, and U.~Hahn.
\newblock Global envelope tests for spatial processes.
\newblock {\em J. R. Stat. Soc. Ser. B. Stat. Methodol.}, 79(2):381--404, 2017.

\bibitem{owada}
T.~Owada and A.~Thomas.
\newblock Limit theorems for process-level {B}etti numbers for sparse,
  critical, and {P}oisson regimes.
\newblock {\em arXiv preprint arXiv:1809.05758}, 2018.

\bibitem{peccati}
G.~Peccati and M.~S. Taqqu.
\newblock {\em Wiener Chaos: Moments, Cumulants and Diagrams}.
\newblock Springer, Milan, 2011.

\bibitem{mc2}
A.~H. Rafati, F.~Safavimanesh, K.-A. Dorph-Petersen, J.~G. Rasmussen,
  J.~M{\o}ller, and J.~R. Nyengaard.
\newblock Detection and spatial characterization of minicolumnarity in the
  human cerebral cortex.
\newblock {\em Journal of Microscopy}, 261(1):115--126, 2016.

\bibitem{gibbsCLT}
A.~Xia and J.~E. Yukich.
\newblock Normal approximation for statistics of {G}ibbsian input in geometric
  probability.
\newblock {\em Adv. in Appl. Probab.}, 47(4):934--972, 2015.

\bibitem{yogeshAdler1}
D.~Yogeshwaran and R.~J. Adler.
\newblock On the topology of random complexes built over stationary point
  processes.
\newblock {\em Ann. Appl. Probab.}, 25(6):3338--3380, 2015.

\bibitem{yogeshAdler2}
D.~Yogeshwaran, E.~Subag, and R.~J. Adler.
\newblock Random geometric complexes in the thermodynamic regime.
\newblock {\em Probab. Theory Related Fields}, 167(1-2):107--142, 2017.

\end{thebibliography}
\bibliographystyle{abbrv}

\section{Proof of Theorem \ref{pdThm}}
\label{proofSec}
The main tool to prove Theorem \ref{pdThm} is the general CLT of \cite[Theorem 1.14]{yogeshCLT}. To be in that framework, we need to express the quantity $\lan f,\PDDn \ran = \sum_{ i \in I^{M,q}(\PP_n)} f\big(B_i^M, D_i^M \big)$ in the form
$\sum_{\x \in \PP_n } \xi(\x, \PP_n)$ for a suitable score function $\xi(x, \PP_n)$.

In other words, we need to transform the indexing over features into an indexing over the points of the point process $\PP_n$. We achieve this goal by assigning to each feature a point $x \in \PP_n$ that either kills or gives birth to this feature, depending on whether $q = 0$ or $ q = 1$. 

First, the death of a cluster at time $r > 0$ is always caused by the merging of two points $x, x' \in \PP_n$ at distance $2r$. Indeed, when the size of a component has a jump, this can only appear by attaching to another component. If $\CO_r(x)$ dies by this merging, we say that $x'$ \emph{kills} $\CO_r(x)$. This ensures that if two components both die when they merge, their deaths are caused by different points.   

Similarly, if $q = 1$, then the birth of a hole at time $r > 0$ is caused by two points $x, x' \in \PP_n$ at distance $2r$ whose connection creates a new hole. If only one feature is born at time $r$, we choose the lexicographic minimum of $x$ and $x'$ and say that it \emph{gives birth} to this hole. However, if a large hole is split into two $M$-bounded holes, it can happen that two holes $H,H'$ are born at the same time. In this situation, we assign one hole to each of $x$ and $x'$. Hence, we define the score functions as 
	\begin{equation}
		\begin{split}
        \label{xiDefEq}
	\xi_0(\x, \PP_n) &= \sum_{ i \in I^{M, 0}(\PP_n)} \one\{\text{$\x$ kills the $i$th cluster}\} f(0, D_i^M),\\
	\xi_1(\x, \PP_n) &= \sum_{ i \in I^{M, 1}(\PP_n)} \one\{\text{$\x$ gives birth to the $i$th hole}\} f(B_i^M, D_i^M).
		\end{split}
	\end{equation}

%
%
Definition \eqref{xiDefEq} translates the desired CLT for $\lan f, \PDDn \ran$ into the framework of \cite[Theorem 1.14]{yogeshCLT}. It remains to verify the conditions stated therein. 
\begin{proof}[Proof of Theorem \ref{pdThm}]
	According to \cite[Theorem 1.14]{yogeshCLT} we have to verify that the pair $(\PP, \xi_q)$ belongs to class (A2) (see Definition \cite[Definition 1.7]{yogeshCLT}) and that the $p$th moment condition \cite[Equation (1.19)]{yogeshCLT} holds for every $p$. 
	
	%
	%
	Belonging to class (A2) involves itself three conditions. The first is exponential decay of correlations, one of our standing assumptions on the point process $\PP$. The second asks for an exponentially decaying radius of stabilization. Since we work with $M$-bounded features, this radius is finite. Finally, we need to verify the power-growth condition \cite[Equation (1.18)]{yogeshCLT} stating that for $\Wr(x) = x + [-\sqrt r/2, \sqrt r/2]^2$, the upper bound
\begin{align*}
	\xi_q(x, \X \cap \Wr(x)) \one\{ \#(\X \cap \Wr(x)) = n\} \le c^n (1 \vee r^n)
\end{align*}
	holds for every $r > 0$, locally finite $\X \subset \R^2$ and $x \in \X$. To achieve this goal, we note that in the worst case $x$ can be responsible for the death of all other points of $\X$. Similarly, it can give birth to at most $\X(\Wr(x)) - 1$ holes. Hence,
	$$\xi_q(\x, \X \cap \Wr(x)) \one\{\#(\X \cap \Wr(x)) = n\} \le |f|_\infty (n - 1) \le (1 + |f|_\infty)^n.$$

%
%
Finally, we verify the $p$th moment condition \cite[Equation (1.19)]{yogeshCLT}. That is, we prove that for every $p > 0$ there exists $M_p > 0$ such that
\begin{align}
        \label{momBoundEq}
	\sup_{\substack{n \ge 1 \\ l \le p,\, \xx \in \R^{2l}}}\E_{\xx}^{!}[|\xi_q(\x_1, \PP_n \cup \xx)|^p] \le M_p.
\end{align}
	We explain in detail how this is achieved if $q = 0$, noting that the case $q = 1$ can be deduced after minor modifications. If $x \in \PP$ is  responsible for the death of a component at time $r$, then there exists $x' \in \PPn$ at distance $2r$ from $x$. Since each ball grows for time at most $\tf$, we see that 
	$$|\xi_0(x, \PPn)| \le |f|_\infty (\PPn(B_{2\tf}(x)) + p)$$
	and an application of condition {\bf (M)} concludes the proof.
       \end{proof}

\section{Proofs of Theorem \ref{mainThm} and Corollary \ref{apfCor}}
\label{verifSec}
%
%
In the following, we assume $q = 1$, since the proofs for $q = 0$ are similar but easier. Hence, to simplify notation, we write $\PD_{b, d}( \PPn)$ for $\PD^{M, 1}_{b, d}( \PPn)$.

\begin{proof}[Proof of Corollary \ref{apfCor}]
	Note that if $(X(s))_{s \le \tf}$ is a Gaussian process, then the process $(\int_0^rX(s)\d s)_{r \le \tf}$ is also Gaussian.
 As mentioned above, the plan is to start from Theorem \ref{mainThm} and then apply the continuous mapping theorem \cite[Theorem 4.27]{kallenberg}. To this end, we show that $\{\APFo_r(\PPn)\}_{r \le \tf}$ is a continuous functional of the persistent Betti numbers $\{\PD_{b, d}(\PPn)\}_{b, d \le \tf}$. We assert that
\begin{align}
\label{apfEq}
	\APFo_r(\PPn) = \int_0^r \PD_{b, 0}(\PPn) \d b + \int_0^{\tf} \PD_{r, t}(\PPn) \d t - r\PD_{r, 0}(\PPn).
\end{align}
	The remainder of the proof proceeds in two steps. First, we verify identity \eqref{apfEq}. Second, we show that the right-hand side is continuous in $\PD$ with respect to the Skorokhod topology.

	To prove identity \eqref{apfEq}, linearity allows us to reduce the claim to the case where the persistence diagram consists of a single $\de$-measure at a point ${(B_0, D_0)}$ for some $D_0 > B_0 > 0$. If $B_0 > r$, then both sides vanish. If $B_0 \le r$, then $\PD_{b, 0} = \one\{b \ge B_0\}$ and $\PD_{r, t} = \one\{t \le D_0\}$, so that the right-hand side of \eqref{apfEq} gives the asserted
\begin{align*}
(r - B_0) + D_0 - r = (D_0 - B_0).
\end{align*}

	Let $\b \in D(\tfs, \R)$, where $D(\tfs, \R)$ is the Skorokhod space of c\`adl\`ag functions from $\tfs$ to $\R$. For any $r \ge 0$ put
	$$\Phi_r(\b) =  \int_0^r \b_{b, 0} \d b + \int_0^{\tf} \b_{r, t} \d t - r\b_{r, 0}.$$
	According to \eqref{apfEq}, it is sufficient to prove that the function $\Phi_r:\,D(\tfs, \R) \to D([0, \tf], \R)$, $\b \mapsto (\Phi_r(\b))_{r \le \tf}$ is continuous with respect to the Skorokhod topology. We prove this for the first integral. The arguments for the second are similar. Let $\PD':\, \tfs \to \R$ be c\`adl\`ag and $\la: \tfo \to \tfo$ be an increasing continuous bijection. Then,
	\begin{align*}
		&\Big| \int_0^{\la(r)} \PD_{b, 0} \d b - \int_0^r \PD'_{b, 0} \d b\Big| \le |\la(r) - r||\PD_{\cdot, 0}|_\infty + \int_0^r|\PD_{b, 0}  - \PD'_{b, 0}| \d b\\ 
		&\quad\le |\la(r) - r||\PD_{\cdot, 0}|_\infty + \int_0^r|\PD_{\la(b), 0}  - \PD'_{b, 0}|\d b + \int_0^r|\PD_{\la(b), 0}  - \PD_{b, 0}|\d b.
	\end{align*}
	If $\PD'$ approaches $\PD$ in the Skorokhod metric, then by definition of this metric, we can choose $\la$ such that the first two expressions become arbitrarily small. Moreover, since $\PD$ itself is c\`adl\`ag, it follows that also the third expression tends to 0.
\end{proof}

%
%
The proof of Theorem \ref{mainThm} decomposes into two steps: lower and upper variance bounds and an upper bound on fourth-order cumulants.  In what follows, we write
$$\PD(E, \PPn) = \PD_{b_+, d_+}( \PPn) + \PD_{b_-, d_-}(\PPn) -\PD_{b_+, d_-}(\PPn) - \PD_{b_-, d_+}( \PPn)$$
for the increment of $\PD_{b, d}$ in the block $E = (b_-, b_+] \times (d_-, d_+]$ with $b_- < b_+$ and $d_- < d_+$. Notice that this is minus the measure $\PDD(\PPn)$ from \eqref{pdEq} evaluated at the block $E$. Moreover, $\PD(E,\PPn)$ is the number of holes with birth time before $b_+$ and death time between $d_-$ and $d_+$ minus the number of holes with birth time before $b_-$ and death time between $d_-$ and $d_+$. Following \cite{bickel}, two blocks $E, E' \subset \tfs$ are \emph{neighboring} if they share a common side.

%
%
\begin{proposition}[Variance lower bound]
\label{varProp1}
	Let $\PP$ be a conditionally $m$-dependent point process with exponential decay of correlations. Moreover, let $a_1,\dots, a_k \ne 0$ and $E_1,\dots, E_k \subset [0, \tf]^2$ be pairwise disjoint blocks such that each $E_i$ contains some $(b, d) \in \tfs$ with $d > b$. 
	Then,
			$$\liminf_{n \to \infty} \frac1n \Var\Big(\sum_{i \le k} a_i \PD(E_i, \PPn)\Big) > 0.$$
\end{proposition}

\begin{proposition}[Variance upper bound]
\label{varProp}
	Let $\PP$ be a conditionally $m$-dependent point process with exponential decay of correlations. Then, there exist $n_0 \ge 1$ and $\e_0, C_0 > 0$ such that
	$$\frac1n \Var\big(\PD(E, \PPn)\big) \le C_0 |E|^{1/2 + \e_0}$$
holds for all $n \ge n_0$ and blocks $E \subset \tfs$.
\end{proposition}

Now, the \emph{$k$th cumulant $c^k$} of $k \ge 1$ real random variables $Y_1, \dots, Y_k$ equals
$$c^k(Y_1, \dots, Y_k) = \sum_{\{T_1, \dots, T_p\} \preceq \{1, \dots, k\}}(-1)^{p - 1} (p - 1)!\,\E\Big[\prod_{i \in T_1} Y_i\Big] \cdots \E\Big[\prod_{i \in T_p} Y_i\Big],$$
provided that all appearing moments are well-defined \cite[Proposition 3.2.1]{peccati}. Here, the sum ranges over all partitions $\{T_1, \dots, T_p\}$ of the set $\{1, \dots, k\}$. 

%
%
\begin{proposition}[Cumulant bound]
\label{tightProp}
	Let $\PP$ be a conditionally $m$-dependent point process  with exponential decay of correlations  satisfying conditions {\bf(AC)} and {\bf (M)}. Then, there exist $n_0' \ge 1$ and $\e_0', C_0' > 0$ such that
    \begin{align*}
\frac1n c^4\big(\PD(E, \PPn), \PD(E, \PPn),\PD(E', \PPn), \PD(E', \PPn)\big) 
		\le C_0' |E|^{1/2 + \e_0'}|E'|^{1/2 + \e_0'}
    \end{align*}
holds for all $n \ge n_0'$ and neighboring blocks $E, E' \subset \tfs$.
\end{proposition}
We postpone the proofs of Propositions \ref{varProp1}--\ref{tightProp} to Sections \ref{varSec}--\ref{tightSec2}, respectively. To deduce Theorem \ref{mainThm} from these two central auxiliary results, we write
$$\mPD_{b, d}( \PPn) = \PD_{b, d}( \PPn) - \E[\PD_{b, d}( \PPn)]$$
for the centered persistent Betti numbers.

%
%
\begin{proof}[Proof of Theorem \ref{mainThm}]
	Let $a_1',\dots, a_{k'}' \ne 0$ and $(b_1, d_1), \dots, (b_{k'}, d_{k'}) \in \tfs$ be pairwise distinct, and put 
	$$X_n = \sum_{i \le k'} a_i' \PD_{b_i', d_i'}( \PP_n).$$
	Then, after suitable regrouping of terms, we can express $X_n$ in the form
	$$X_n = \sum_{i \le k} a_i \PD(E_i, \PPn).$$
as in Proposition \ref{varProp1}. Now, combining Proposition \ref{varProp} with Theorem \ref{pdThm} and the variance asymptotics \cite[Theorem 1.12]{yogeshCLT} shows that the centered and rescaled random variable $n^{-1/2}(X_n - \E[X_n])$ converges in distribution to a Gaussian. Hence, the Cram\'er-Wold device yields convergence of the finite-dimensional distributions of $n^{-1/2}\mPD_{b, d}( \PP_n)$. 

	Next, \cite[Proposition 3.2.1]{peccati} gives the general cumulant identity 
	\begin{align*}
		\E[X^2 Y^2] &= c^4(X, X, Y, Y) + \Var(X)\Var(Y) + 2 \ms{Cov}(X, Y)^2 \\
		&\le c^4(X, X, Y, Y) + 3\Var(X)\Var(Y) 
	\end{align*}
	for centered random variables $X, Y$. Hence, by Propositions \ref{varProp} and \ref{tightProp},
	\begin{align*}
		\E\big[n^{-2}\mPD(E, \PPn)^2\mPD(E', \PPn)^2 \big] \le (C_0'/n + 3C_0^2)|E|^{1/2 + \e_0''}|E'|^{1/2 + \e_0''},
	\end{align*}
	for some $\e_0'' > 0$.
	In particular, the process $\big\{n^{-1/2}\mPD_{b, d}( \PP_n)\big\}_{b, d \le \tf}$ is tight in Skorokhod topology \cite[Lemma 3]{heinrichShot}. In this context, we note that condition (8.4) of \cite[Lemma 3]{heinrichShot} follows from the variance upper bound derived in Proposition \ref{varProp} and that similar as in (2.18) \cite{kclt}, we have replaced the equality in (8.5) of \cite[Lemma 3]{heinrichShot} by an inequality.
	Combining this property with the convergence of finite-dimensional distributions yields the asserted weak convergence.
\end{proof}

\subsection{Proof of Proposition \ref{varProp1}}
\label{varSec}

To show the variance lower bound, we adapt a conditioning argument that has already been successfully applied in the setting of Gibbsian point processes \cite{gibbsCLT}. More precisely, we subdivide the window $W_n$ into blocks of a fixed size and use the law of conditional variance to obtain a lower bound in the order of the number of blocks.

Associate with the $j$th feature $H_j$ in $\ms{PD}^{M, 1}(\PP_n)$ a center point $y_j \in W_n$, for instance by taking the point $p(H_j)$ as defined in Section \ref{loopSec}. Then, 
$$\nun = \sum_{i \le k}a_i \sum_{j \in I^{M, 1}(\PP_n)} \one\{(B_j^M, D_j^M) \in E_i\} \de_{y_j}$$
defines a signed measure of total mass $\nun(\Wn) = \sum_{i \le k} a_i \PD(E_i, \PPn)$.

%
%
In the vein of \cite{gibbsCLT}, the key towards proving a lower bound on the variance is the following non-degeneracy property, where $\rac$ is introduced in Section \ref{resultSec}.
\begin{lemma}[Non-degeneracy]
	\label{nDegLem}
	It holds that 
	$$\inf_{\substack{n \ge t \ge \rac^2}}\E\big[\Var\big(\nun(W_t)|\, \La, \PP \setminus W_{\rac^2}\big)\big] > 0.$$
\end{lemma}

Before proving Lemma \ref{nDegLem}, we explain how it implies Proposition \ref{varProp1}. In essence, the proof follows along the lines of \cite[Lemma 4.3]{gibbsCLT}. Nevertheless, since the details of the conditioning argument differ a bit from the corresponding picture for Gibbs processes, we explain how to adapt the main steps from \cite[Lemma 4.3]{gibbsCLT} in the present setting.
\begin{proof}[Proof of Proposition \ref{varProp1}]
    The idea of proof is to consider a family of well-separated blocks
	in $\Wn$. Then, we leverage the conditional $m$-dependence of the point process and the $M$-boundedness of the features to decompose the variance of their contributions as the sum of the variances. More precisely, we apply the assumption of conditional $m$-dependence with the conditioning set
	$$A'' = \R^2 \setminus \bigcup_{z\in \Z^2} (6\rho z + W_{\rho^2})$$ 
	chosen as the complement of the union of well-separated blocks of side length $\rho = m \vee \rac$. Then, the law of total variance yields the lower bound 
	$$\Var(\nun(\Wn)) \ge \E\big[\Var\big(\nun(\Wn)\, |\, \La, \PP \cap {A''}\big)\big].$$
	Moreover, since $\rho > M$ the statistics $\nun((A'')^-)$ in the smaller domain 
	$$(A'')^- = \R^2 \setminus \bigcup_{z\in \Z^2} (6\rho z + W_{9\rho^2})$$ 
	is measurable with respect to $\PP \cap {A''}$. We obtain that 
		$$\E\big[\Var\big(\nun(\Wn) | \La, \PP \cap {A''}\big)\big]  = \E\big[\Var\big(\nun(\R^2 \setminus (A'')^-) | \La, \PP \cap {A''}\big)\big]$$
		because $\nun((A'')^-)$ is $\PP \cap A''$ measurable. Thanks to the conditional $m$-dependence, we have 
	\begin{align*}\E\big[\Var\big(\nun(\Wn) | \La, \PP \cap {A''}\big)\big]  &= \sum_{\substack{z \in \Z^2 \\ 6\rho z + W_{9\rho^2} \subset \Wn}} \E\big[\Var\big(\nun(6\rho z + W_{9\rho^2})\, |\, \La, \PP \cap {A''}\big)\big]\\
		&\ge \hspace{-.8cm} \sum_{\substack{z \in \Z^2 \\ 6\rho z + W_{9\rho^2} \subset \Wn}}\hspace{-.7cm} \E\big[\Var\big(\nun(6\rho z + W_{9\rho^2})\, |\, \La, \PP \setminus (6\rho z + W_{\rho^2})\big)\big].
	\end{align*}
	Now, the number of $6\rho$-blocks contained in $\Wn$ is of order $n$, and we conclude by noting that Lemma \ref{nDegLem} and $\rho > \rac$ imply that each of the contributions is bounded away from 0.
\end{proof}

To verify non-degeneracy, we rely on the techniques introduced in \cite{gibbsCLT}. In particular, we make use of \cite[Lemma 2.3]{gibbsCLT}, which we restate below to render the presentation self-contained.
%
%
\begin{lemma}
	\label{23Lem}
	Let $Y$ be a real random variable and $A_1, A_2$ be Borel sets of $\R$. Then, 
	$$\Var(Y) \ge \frac14 \min_{i \in \{1, 2\}}\P(Y \in A_i)\inf_{x_1 \in A_1, x_2 \in A_2}|x_1 - x_2|^2 .$$
\end{lemma}

\begin{proof}[Proof of Lemma \ref{nDegLem}]
	Write 
	$$F_1 = \{\PP \cap W_{\rac^2/9} = \es\}\quad\text{and}\quad F' = \{\PP \cap (W_{\rac^2} \setminus W_{\rac^2/9}) = \es\}$$
	for the events that there are no points in $W_{\rac^2/9}$ and $W_{\rac^2} \setminus W_{\rac^2/9}$, respectively. Next, let 
	$$F_2 = \{\PD(E_1, \PP_{\rac^2/9}) = 1\} \cap \{\PD(E_2, \PP_{\rac^2/9}) = \cdots = \PD(E_k, \PP_{\rac^2/9}) = 0\}$$
	denote the event 
	that all but the first of the considered persistent Betti numbers vanish. 
	Now, let $I_0$ denote the indices of all features that are entirely contained in $\R^2 \setminus W_{\rac^2}$ and put 
	$$Y = \sum_{i \le k}a_i \#\{{j \in I^{M, 1}(\PP_n) \setminus I_0(\PP_n)}:\, (B_j^M, D_j^M) \in E_i\}.$$
	Then, by Lemma \ref{23Lem} with $A_1 = [a_1, \infty)$ and $A_2 = \{0\}$,
	\begin{align*}
		\E\big[\Var\big(\nun(W_t)|\, \La, \PP \setminus W_{\rac^2}\big)\big] &= \E\big[\Var\big(Y|\, \La, \PP \setminus W_{\rac^2}\big)\big] \\  &\ge \frac{a_1^2}4\E\big[\min_{i \in \{1, 2\}}\P(\PP \in F' \cap F_i | \La, \PP \setminus W_{\rac^2})\big],
	\end{align*}
	and it remains to show that the right-hand side is non-zero.

Since $E_1, \dots, E_k$ are pairwise disjoint and contain points above the diagonal, \cite[Example 1.8]{shirai} shows that under the homogeneous Poisson point process the event $F' \cap F_2$ has positive probability. Also $F' \cap F_1$ is of positive probability. Hence, an application of condition {\bf (AC)} concludes the proof.
\end{proof}

\subsection{Proof of Proposition \ref{varProp}}
\label{tightSec}


%
%
For a block $E = (b_-, b_+] \times (d_-, d_+] \subset \tfs$, we let $\xi_E$ denote the score function associated with $\PD(E, \PPn)$. That is, 
$$\xi_E(x, \PP_n) = \#\{(\bim, \dim) \in E:\, \text{$x$ gives birth to the $i$th hole}\}$$
is the number of holes born by $x$ with birth and death times in $E$. Note that if $x$ gives birth to the $i$th hole, then it gets in contact with another point at time $\bim \in (b_-, b_+]$. In particular, $\PP$ contains a point in the annulus $A_{2b_-, 2b_+}(x) = B_{2b_+}(x) \setminus B_{2b_-}(x)$.

 Moreover, if the $i$th hole dies at time $\dim \in (d_-, d_+]$, then a previously vacant component is covered completely, which is caused by three disks centered at points in $\PP$ meeting at a single point in the plane. The three center points of the disks must form a triangle with no obtuse angle. Otherwise, two of the disks would meet for the first time in the interior of the third and hence no connected component in the background was covered by the merging. This could be interpreted as a feature that is born and dies at the same time, but we chose to exclude such features in our definition of 1-features.
 
 Henceforth, let $B^\pm_d(x, y) \subset \R^2$ denote the two disks of radius $d > 0$ whose boundary passes through $x, y \in \R^2$. If $|x - y|/2 > d$, we let $B^\pm_d(x, y)$ be empty. The points in $B^+_d(x, y) \cup B^-_d(x, y)$ are exactly the points $z$ such that the time when the boundaries of the three disks around $x$, $y$, and $z$ meet in one point is at most $d$. 
 For $d_+ > d_-\ge 0$, we let
\begin{equation*}
 D_{d_-, d_+}(x, y) = \big( B^+_{d_+}(x, y) \cup B^-_{d_+}(x, y) \big) \backslash \big( B^+_{d_-\vee a}(x, y) \cup B^-_{d_-\vee a}(x, y)\big),
\end{equation*}
where $a = |x-y|/2$.
This set consists of all points $z$ such that the boundaries of the three disks around $x$, $y$ and $z$ meet at time $r$ with $d_- < r \le d_+$. Some $z\in D_{d_-, d_+}(x, y)$ may still form a triangle having an obtuse angle with $x$ and $y$, that is, the disks around $x$, $y$, and $z$ already met earlier in an interior point of one of the disks. However, all $z$ that can cause the death of a hole in $E$ together with $x$ and $y$ must be contained in $D_{d_-, d_+}(x, y)$.

Now,
\begin{align}
	\label{annBoundEq}
	\xi_E(x, \PP_n) \le \PP(B_M(x)) \one_{E_x},
\end{align}
where $E_x$ denotes the event that for some $\PP' \subset \PP$ with $x \in \PP'$ the event $E_{x, \ms b}(\PP')\cap E_{x, \ms d}(\PP')$ occurs, where 
\begin{align*}
	E_{x, \ms b}(\PP') ={}& \big\{x \text{ creates an $M$-bounded hole in $\{U_r(\PP')\}_{r\geq 0}$
		with birth and}\\ &\text{\quad death time in $E$ by connecting to some $x_1 \in \PP' \cap A_{2b_-, 2b_+}(x)$}\big\}\\
	E_{x, \ms d}(\PP') ={}&\big\{\text{$\exists y_1, y_2 \in \PP' \cap B_M(x)$ such that $y_1, y_2, y_3$ kill an $M$-bounded hole }\\
		&\text{\quad in $\{U_r(\PP')\}_{r\geq 0} $ with birth and death time in $E$ for some}\\
		&\text{\quad $y_3 \in \PP' \cap D_{d_-, d_+}(y_1, y_2)$}\big\}.
\end{align*}
Here, we say that $y_1, y_2, y_3\in \PP'$ kill the hole $H$ if the disks around the points meet for the first time at $p(H)$. In particular, any three points can kill at most one hole.

Similarly, for a block $E' = (b_-', b_+'] \times (d_-', d_+'] \subset \tfs$, 
\begin{align*}
	\xi_E(x, \PP_n) \xi_{E'}(x', \PP_n)\le \PP(B_M(x))\PP(B_M(x')) \one_{E_{x, x'}''},
\end{align*}
where we let $E_{x, x'}''$ denote the event that for some $\PP' \subset \PP$ with $x, x' \in \PP'$ the event 
$$E_{x, \ms b}(\PP')\cap E_{x', \ms b}'(\PP')\cap E_{x, \ms d}(\PP') \cap E_{x', \ms d}'(\PP')$$
occurs.

%
%
Using this notation, the proof of the variance upper bound is now based on the following pivotal geometric moment bound. In the following, $\P_{\xx}$ denotes the {unreduced Palm measure} characterized via
$\E_{\xx}[f(\PP)] = \E_{\xx}^![f(\PP \cup \xx)]$ for any non-negative measurable $f:\, \NN \to [0, \infty)$. We recall from \eqref{factMomEq} that $\rho^{(p)}$ denotes the $p$th factorial moment density. In the following, we adhere to the convention $\int_{B^0} f(x) \d z = f(x)$.

\begin{lemma}[Moment bound]
    \label{geom1Lem}
    Let $\PP$ be a stationary point process having fast decay of correlations and satisfying condition {\bf (M)}.
    Let $p \ge 0$ and $K_0 > 0$. Then, there exist $\e>0$ and $C_{\ms g} > 0$ such that for all $n>0$ and any ball $B \subset \R^2$ of radius $K>K_0$,
    \begin{enumerate}
        \item \label{moment1}
        $$
            \frac 1{|B|^{p + 1}}\int_{B^p} \P_{o, \zz}(E_o) \rho^{(p + 1)}(o, \zz) \d \zz \le C_{\ms g} |E|^{1/2 + \e}
        $$
        holds for all blocks $E \subset \tfs$.
        \item \label{moment2}
		$$\frac 1{|B|^{p + 1 }}\int_{B^{p }} \P_{o, \zz}(E_{o, o}'' ) \rho^{(p + 1)}(o, \zz) \d \zz \le C_{\ms g}|E|^{1/2 + \e}|E'|^{1/2 + \e}$$
        holds for all neighboring blocks $E, E' \subset \tfs$, and
        \item \label{moment3}
            $$\frac 1{|B|^{p+2 }}\int_{B^{p + 1 }} \P_{o, z', \zz}(E''_{o, z'}) \rho^{(p + 2)}(o, z', \zz) \d (z', \zz) \le C_{\ms g}|E|^{1/2 + \e}|E'|^{1/2 + \e}$$
            holds for all neighboring blocks $E, E' \subset \tfs$.
    \end{enumerate}
\end{lemma}

The proof of Lemma \ref{uDiagConLem} relies on a delicate geometric analysis that we defer to Section \ref{geomProofSec}.
We now prove Proposition \ref{varProp}. As in \cite[Equation (1.6)]{yogeshCLT}, for $\xx = (x_1, \dots, x_p) \in \R^{2p}$ and $k_1, \dots, k_p \ge 0$, we introduce the \emph{mixed $\xi_E$-moments}
\begin{align}
	\label{mxmEq}
	m_n^{(k_1, \dots, k_p)}(\xx) = \E_{\xx}[\xi_E(x_1, \PP_n)^{k_1} \cdots \xi_E(x_p, \PP_n)^{k_p}] \rho^{(p)}(\xx).
\end{align}
In the rest of the manuscript, we freely use that exponential decay of correlations implies boundedness of the factorial moment densities \cite[Inequality (1.11)]{yogeshCLT}.

%
%
\begin{proof}[Proof of Proposition \ref{varProp}]
        To lighten notation, we write $\xi$ instead of $\xi_E$. To give the paper more pleasant to read, we have not attempted to optimize the exponents occurring in the course of this proof.
	Proceeding as in \cite[Equation (4.1)]{yogeshCLT}, the refined Campbell-Mecke formula \cite[Equation (1.9)]{yogeshCLT} gives that $\Var\big(\PD(E, \PPn)\big)$ equals
        \begin{align}
                \label{varEq}
                \int_{\Wn} m_n^{(2)}(x)\d x
		+ \int_{\Wn \times \Wn} \big(m_n^{(1,1)}(x,y) - m_n^{(1)}(x)m_n^{(1)}(y)\big)\d(x, y).
        \end{align}
        We derive bounds for the two summands separately.

        %
        %
        By stationarity, \eqref{annBoundEq} and H\"older's inequality, the first expression is at most
        \begin{align}
                \label{m2Eq}
                n \E_o[\PP(B_M(o))^{16}]^{1/8} \P_o(E_o)^{7/8} \rho = n \big(\E_o[\PP(B_M(o))^{16}]\rho\big)^{1/8} \big(\P_o(E_o)\rho\big)^{7/8}.
        \end{align}
        Hence, Lemma \ref{geom1Lem}(1) with $p = 0$ yields the asserted upper bound.

        %
        %
        To deal with the double integral in \eqref{varEq}, we recall that $\xi$ is a local score function and that $\PP$ exhibits exponential decay of correlations. Hence, as in \cite[Equation (3.26)]{yogeshCLT}, the factorial moment measure expansion shows that
        $$\big|m_n^{(1,1)}(x,y) - m_n^{(1)}(x) m_n^{(1)}(y)\big| \le c\phi(|x - y|/2)$$
	for some $c > 0$. In particular, choosing a cut-off $K = |E|^{-1/128}$, we see that 
	$$\sup_{x \in \Wn}\int_{\Wn \setminus B_K(x)} \big|m_n^{(1,1)}(x,y) - m_n^{(1)}(x)m_n^{(1)}(y)\big|\d y \le C |E|$$
        holds for a suitable $C > 0$ and it suffices to derive an upper bound for
        \begin{align*}
                \int_{B_K(x)}m_n^{(1,1)}(x,y) + m_n^{(1)}(x) m_n^{(1)}(y)\d y.
        \end{align*}
	For the second summand, we can argue similarly as in \eqref{m2Eq}, so that it remains to bound the integral involving $m_n^{(1,1)}(x,y)$. Here, we set $z = y - x$, note that $\rho^{(2)}(x, y) = \rho^{(2)}(o, z)$ and combine \eqref{annBoundEq} with H\"older's inequality to arrive at
        $$m_n^{(1, 1)}(x, y) \le \E_{o, z}[\PP(B_M(o))^{16}]^{1/16}\E_{o, z}[\PP(B_M(z))^{16}]^{1/16} \P_{o, z}(E_o)^{7/8} \rho^{(2)}(o, z).$$
	We bound $\E_{o, z}[\PP(B_M(z))^{16}]$ thanks to condition {\bf (M)}.
	Finally, by Jensen's inequality applied to the uniform distribution on $B_K(o)$,
	\begin{align}
		\label{diagEq}
		\frac 1{|B_K(o)|}\int_{B_K(o)}\hspace{-.5cm}(\P_{o, z}(E_o)\rho^{(2)}(o, z))^{7/8}\d z \le\Big(\frac 1{|B_K(o)|}\int_{B_K(o)}\hspace{-.4cm}\P_{o, z}(E_o)\rho^{(2)}(o, z)\d z\Big)^{7/8}\hspace{-.4cm},
	\end{align}
	so that applying Lemma \ref{geom1Lem}(1) with $p = 1$ shows that the right-hand side is of order at most $(|E|^{3/4}|B_K(o)|)^{7/8} = |E|^{7 \cdot 47 / 512}$, thereby concluding the proof.
\end{proof}

\subsection{Proof of Proposition \ref{tightProp}}
\label{tightSec2}

%
%

To prove Proposition \ref{tightProp}, we take up the idea suggested in \cite[Theorem 8]{heinrichShot} and \cite[Theorem 8.1]{calka} and express $c^4$ in terms of cumulant measures induced by the functional of interest. A slight technical nuisance in the present setting comes from dealing with a product of two different functionals -- one associated with the block $E$ and the other with $E'$ -- whereas the semi-cluster measure machinery from \cite[Section 4.3]{yogeshCLT} relies on a single score function. However, this artificial difficulty can be overcome by formally attaching $\{1, 2\}$-valued marks to $\PPn$. Taking up the notation from \cite{raic3}, we let $\bR^2 = \R^2 \times \{1, 2\}$ and $\bPPn = \PPn \times \{1, 2\}$ denote the correspondingly marked space and point process.
Writing $E'' = (E, E')$, we define an augmented score function $\xi_{E''}$, where points with mark 1 are evaluated with the first score function and points with mark 2 are evaluated with the second score function. In other words, 
$$\xi_{\mm}((x, \tau), \bPPn)=
\begin{cases}
	\xi_E(x, \PPn)&\text{ if $\tau = 1$, }	\\
	\xi_{E'}(x, \PPn)&\text{ if $\tau = 2$.}	
\end{cases}$$

We take the concise proof of Proposition \ref{varProp} as a blueprint for the strategy of the more involved setting laid out in Proposition \ref{tightProp}. In particular, we need to address two main steps: bounds for mixed moments and a reduction of the integral to the diagonal. 

%
%
In order to reduce to the diagonal, we decompose the cumulant measure into semi-cluster measures as in \cite[Section 5.1]{bary} and \cite[Section 3.2]{raic3}. For the convenience of the reader, we reproduce the basic definitions. First, the \emph{$k$th moment measure $M^k(\mun)$} is given as
$$\lan \ff, M^k(\mun) \ran = \int \ff(\bxx) M^k(\mun)(\d \bxx)  = \E[\lan f_1, \mun\ran \cdots \lan f_k, \mun \ran], $$
where $\ff = f_1 \otimes \cdots \otimes f_k$ is non-negative and measurable with each $f_i$ defined on $\bR^2$, and
$$\mun = \mu_{\mm, n} = n^{-1} \sum_{\bx \in \bPPn}\xi_{\mm}(\bx, \bPPn) \de_{\bx}$$
denotes the empirical measure associated with $\xi_{\mm}$ and $\bPPn$.
In terms of mixed $\xi$-moments, with $\bxx_{T_i}$ the projection of $\bxx$ to the coordinates in $T_i$, we write
\begin{align}
	\label{singExpEq}
	\d M^k= \sum_{\{T_1, \dots, T_p\} \preceq \{1, \dots, k\}}m_n^{(T_1, \dots, T_p)} \d \bxx_{T_1} \cdots \d \bxx_{T_p}, 
\end{align}
	where $\d \bxx_{T_i}$ are the singular differentials determined via 
	$$\int_{\bR^{2|T_i|}} f(\bxx_{T_i}) \d \bxx_{T_i} = \int_{\bR^2} f(\bx, \dots, \bx)\d \bx$$ 
	where $f:\bR^{2|T_i|} \to [0, \infty)$ is any non-negative measurable function \cite[Section 3.1]{raic3}. As in \eqref{mxmEq}, for  $k_1, \dots, k_p \ge 0$, the {mixed $\xi_{E''}$-moments} are given as
\begin{align*}
	m^{(T_1, \dots, T_p)}_n(\bxx) = \E_{\xx}[\xi_{E''}(\bx_1, \PP_n)^{|T_1|} \cdots \xi_{E''}(\bx_k, \PP_n)^{|T_p|}] \rho^{(p)}(\xx),
\end{align*}
for every $\bxx = ((x_1, \tau_1), \dots, (x_k, \tau_k)) \in \bR^{2k}$.
	
	Similarly, the \emph{$k$th cumulant measure $\ckn = c^k(\mun)$} equals
$$\lan \ff, \ckn \ran = c^k(\lan f_1, \mun\ran, \dots, \lan f_k, \mun \ran), $$
 so that
\begin{align}
	\label{momExpEq}
	\ckn = \sum_{\{T_1, \dots, T_p\} \preceq \kk}(-1)^{p - 1} (p - 1)!\, M^{T_1} \cdots M^{T_p}, 
\end{align}
where 
	$$\d M^{T_i} = \sum_{\{T_1', \dots, T_{p'}'\} \preceq T_i }m_n^{(T'_1, \dots, T_{p'}')} \d \bxx_{T_1'} \cdots \d \bxx_{T_{p'}'}$$
 denotes the moment measure with coordinates in $T_i$.

%
%
Next, the space $\bWn^4$ decomposes into a union of subsets according to which coordinate is most distant from the diagonal \cite[Lemma 3.1]{raic3}. More precisely, write
$$D(\bxx) = \max_{\{S, T\} \preceq \four}\ms{dist}(\bxx_S, \bxx_T)$$
for the maximal separation of $\bxx_S$ and $\bxx_T$, where $\ms{dist}(\bxx_S, \bxx_T) = \ms{dist}(\xx_S, \xx_T)$. Then, put
$$\sigma(S, T) = \big\{\bxx = (\bxx_S, \bxx_T) \in \bWn^4:\, D(\bxx) = \ms{dist}(\bxx_S, \bxx_T)\big\} \setminus \Delta.$$
Here, the marks are ignored for the diagonal $\Delta \subset \bWn^4$.
We also put $\Wmn = (\Wn \times \{1\})^2 \times (\Wn \times \{2\})^2$. 
%
%
\begin{lemma}[Off-diagonal bounds]
	\label{uDiagConLem}
	 Let $S, T$ denote a non-trivial partition of $\four$. Then, there exist $n_{S, T} \ge 1$ and $\e_{S, T}, C_{S, T} > 0$ such that
	$$\frac1n \big|c^4_n(\S(S, T) \cap \Wmn)\big| \le C_{S, T} |E|^{1/2 + \e_{S, T}}|E'|^{1/2 + \e_{S, T}}$$ 
	holds for all $n \ge n_{S, T}$ and neighboring blocks $E, E' \subset \tfs$.
\end{lemma}

Before proving Lemma \ref{uDiagConLem}, we elucidate how to deduce Proposition \ref{tightProp}.
%
%
\begin{proof}[Proof of Proposition \ref{tightProp}]
	First, 	integration over the cumulant measure decomposes into a diagonal and an off-diagonal part \cite[Equation (3.28)]{raic3}. That is, 
\begin{align*}
	\frac 1n \lan \ff, \cfn \ran &= \frac 1n\int_\Delta \ff \d \cfn + \frac 1n\sum_{S, T} \int_{\sigma(S, T)} \ff \d \cfn.
\end{align*}
	where $\ff = \one_{\Wmn}$ is the indicator function of the domain $\Wmn$ and the sum is over all non-trivial partitions $S, T$. By Lemma \ref{uDiagConLem}, the off-diagonal contributions in this decomposition are bounded above by $\sum_{S, T}C_{S, T}|E|^{1/2 + \e_{S, T}}|E'|^{1/2 + \e_{S, T}}$. 
	
Next, when integrating over the diagonal, we leverage that in the decomposition \eqref{momExpEq}, only $p=1$ contributes \cite[Lemma 3.1]{raic3}. Hence, 
	\begin{align*}
		\int_\Delta \ff \d \cfn &= \int_{W_n }\E_x[\xi_E(\x, \PPn)^2\xi_{E'}(\x, \PPn)^2] \rho \d \x \\
		&\le n\E_o[\PP(B_M)^{2/\e}]^\e \P_o(E_{o, o}'')^{1 - \e}, 
	\end{align*}
	so that applying Lemma \ref{geom1Lem}(2) with $p = 1$ and noting the convention preceding that result concludes the proof.
	\end{proof}

%
%
To prove Lemma \ref{uDiagConLem}, we decompose the cumulant measures into \emph{semi-cluster measures} \cite[Lemma 5.1]{bary}. More precisely, as in \cite{bary, raic3}, any two disjoint non-empty subsets $S', T' \preceq \four$, induce a \emph{cluster measure}
$$U^{S', T'}(A \times B) = M^{S' \cup T'}(A \times B) - M^{S'}(A)M^{T'}(B).$$
Now, $\cfn$ decomposes into semi-cluster measures 
\begin{align}
	\label{semClustEq}
	\cfn = \sum_{\{S', T', T_1, \dots, T_p\} \preceq \four} U^{S', T'}M^{T_1}\cdots M^{T_p}, 
\end{align}
where the sum runs over all partitions such that $S'$ and $T'$ are non-empty subsets of $S$ and $T$, respectively \cite[Lemma 3.2]{raic3}.

Equipped with these ingredients, we now prove Lemma \ref{uDiagConLem}. Since the basic structure of the proof parallels that of Proposition \ref{varProp}, we only provide details for the steps that are substantially different.

%
%
\begin{proof}[Proof of Lemma \ref{uDiagConLem}]
Putting $D_K = \{\bxx \in \Wmn \cap \S(S, T):\, D(\bxx) > K\}$ for $K \ge 1$, we first derive an upper bound for 
	$$\Big|\int_{D_K} \d U^{S', T'}\d M^{T_1} \cdots \d M^{T_p}\Big| = \Big|\int_{D_K} (\d M^{S' \cup T'} - \d M^{S'}\d M^{T'})\d M^{T_1} \cdots \d M^{T_p}\Big|.$$
	For this purpose, we decompose the moment measures $\d M^{S' \cup T'}$, $\d M^{S'}$ and $\d M^{T'}$ according to \eqref{singExpEq}. Hence, we need bounds for the absolute value of differences of mixed $\xi$-moments of the form 
	\begin{equation}
		\label{mixedMBound}
	\begin{aligned}
		&\Big|m_n^{(S''_1, \dots, S''_{p''}, T''_1, \dots, T''_{r''})}(\bxx_{S_1''}, \dots, \bxx_{S_{p''}''}, \bxx_{T_1''}, \dots, \bxx_{T_{r''}''})\\
		&\quad -m_n^{(S''_1, \dots, S''_{p''})}(\bxx_{S_1''}, \dots, \bxx_{S_{p''}''}) m_n^{(T''_1, \dots, T''_{r''})}(\bxx_{T_1''}, \dots, \bxx_{T_{r''}''})\Big|, 
	\end{aligned}
	\end{equation}
	where $\{S''_1, \dots, S''_{p''}\}$ and $\{T''_1, \dots, T''_{r''}\}$ are partitions of $S'$ and $T'$, respectively. Since we are working on the set $\S(S, T)$, as in the proof of Proposition \ref{varProp}, the fast decay of $\xi$-correlations bounds \eqref{mixedMBound} by $c\phi(D(\bxx_{S' \cup T'})/2)$ for a suitable $c > 0$. 
	
	Next, as in \cite[Section 3.1]{raic3} the singular differentials occurring in the expansion \eqref{singExpEq} of the moment measure $M^k$ can be grouped into a single object. More precisely, we write $\tilde \d \bxx$ for the measure that equals $\d \bxx_{T_1} \cdots \d \bxx_{T_p}$ on the subset of $\bR^{2k}$ consisting of all $\bxx = (\bx_1, \dots, \bx_k)$ such that $\bx_i = \bx_j$ if $i, j \in T_r$ for some $r \le p$ and $\bx_i \ne \bx_j$ otherwise. 

	In the setting of the present proof, we note that the bounds on the mixed moments from \eqref{mixedMBound} only involve coordinates with indices in the set $S' \cup T'$. Hence, we need to consider also singular differentials only with respect to these coordinates, i.e., integrate with respect to $\tilde \d \bxx_{S' \cup T'}$. In particular, we arrive at the bound
	\begin{align}
		\label{stEq}
		\Big|\int_{D_K} \hspace{-.4cm} \d U^{S', T'}\d M^{T_1} \cdots \d M^{T_p}\Big| \le c\hspace{-.1cm}\int_{D_K}\hspace{-.4cm} \phi(D(\bxx_{S' \cup T'})/2) \tilde \d \bxx_{S' \cup T'}\d M^{T_1} \cdots \d M^{T_p}.
	\end{align}
	Now, setting $K = |E|^{-\e/128}|E'|^{-\e/128}$, the exponential decay assumption on the function $\phi$ gives control on one integral over the window, while the integrals with respect to the remaining variables are controlled by the volume of balls. Then,	a repeated application of H\"older's inequality provides suitable bounds on the moment measures such that 
	$$\frac1n \int_{D_K} \phi(D(\bxx_{S' \cup T'})/2) \tilde \d \bxx_{S' \cup T'} \d M^{T_1} \cdots \d M^{T_p} \le C|E||E'|$$
	holds for some $C > 0$. Hence, it suffices to provide upper bounds for 
		$$ \frac1n \int_{\{\bxx \in \Wmn:\, D(\bxx) \le K\}} \d M^{T'_1} \cdots \d M^{T'_{p'}}, $$
		where $\{T'_1, \dots, T'_{p'}\}$ is an arbitrary partition of $\four$. We explain how to proceed for $p' = 1$, noting that for $p' > 1$ the arguments are similar but easier.

		We claim that for some $C' > 0$,
			\begin{align}
				\label{uDiagClaimEq}
				\frac1n \int_{W_n \times \{1\}}\int_{B_K(x_1) \times \{1\}} \int_{(B_K(x_1) \times \{2\})^2}\hspace{-.2cm} \d M^{\four} \le C'|E|^{1 + \e/8}|E'|^{1 + \e/8}.
			\end{align}
			To prove this claim, decompose $M^{\four}$ according to \eqref{singExpEq} and let $\{T''_1, \dots, T''_{p''}\}$ be an arbitrary partition of $\four$.
	As in the proof of Proposition \ref{varProp}, a repeated use of H\"older's inequality shows that on $\Wmn$, the mixed moments of the form 
	$$m_n^{(T''_1, \dots, T''_{p''})}(\bx_1, \dots, \bx_4)$$
	are bounded above by $c' \big(\P_{\xx}(E_{x_1, x_i}'')\rho^{(p'')}(\xx)\big)^{1 - \e}$ for a suitable $c' > 0$ and some $i \le 4$. At this point, we may proceed similarly as in \eqref{diagEq} by invoking Lemmas \ref{geom1Lem}(2) and \ref{geom1Lem}(3). As an illustration consider the setting where $p'' = 4$ and $i = 2$. Then, we set $z' = x_2 - x_1$, $z_3 = x_3 - x_1$ and $z_4 = x_4 - x_1$. We combine Jensen's inequality with Lemma \ref{geom1Lem}(3) to show that 
	\begin{align*}
		&\frac1{|B_K|^3}\int_{B_K^3} \big(\P_{o, z', z_3, z_4}(E_{o, z'}'')\rho^{(4)}(o, z', z_3, z_4)\big)^{1 - \e} \d z' \d z_3 \d z_4 \\
		&\quad\le \Big(\int_{B_K^3}\frac1{|B_K|^3}\P_{o, z', z_3, z_4}(E_{o, z'}'')\rho^{(4)}(o, z', z_3, z_4) \d z' \d z_3 \d z_4\Big)^{1 - \e} \\
		&\quad\le C_{\ms g}^{1 - \e}|B_K| |E|^{1/2 + \e/4}|E'|^{1/2 + \e/4}.
	\end{align*}
	Hence, inserting the definition of $K$ concludes the proof.
\end{proof}

\subsection{Proof of Lemma \ref{geom1Lem}}
\label{geomProofSec}

We now turn to the proof of Lemma \ref{geom1Lem}. The proof is based on the following four lemmas that are used to bound the probability with which certain point configurations occur. Throughout we use the notation 
\begin{align*}
E {}&= (b_-, b_+]\times(d_-, d_+],\\
E'{}& = (b_-', b_+']\times (d_-', d_+'],\\
\delta_b{}&=b_+-b_-, \quad \delta_d=d_+-d_-,\\
 \delta_{b'}{}&=b_+'-b_-', \quad \delta_{d'}=d_+' - d_-'.
\end{align*}
	The proofs make use of the inequalities 
\begin{align}\label{sqrtIneq}
|\sqrt x - \sqrt y|{}&\le \sqrt{|x-y|}\\ \label{arcsinIneq}
|\arcsin(x) - \arcsin(y)|{}&\le C_0 \sqrt{|x-y|},
\end{align}
where $C_0>0$ is some constant. Moreover, we repeatedly use that the volume of an annulus is given by
\begin{equation*}
|A_{b_-,b_+}(o)| = b_+^2 - b_-^2 \leq 2b_+\delta_b.
\end{equation*} 


\begin{lemma} \label{AreaDLem}
	Let $x,y\in \R^2$ and $a=|x-y|/2$.  There is a constant $C>0$ such that for all $0\le a \le d_+\le \tf$,
	\begin{align*}
	|D_{d_-,d_+}(x,y)| {}&= 2d_+^2\Big(\pi -  \arcsin\big(\tfrac a{d_+}\big) + \tfrac a{d_+}\sqrt{1-\big(\tfrac a{d_+}\big)^2}\Big) \\
	&-2(d_-\vee a)^2\Big(\pi -  \arcsin\big(\tfrac a{d_-\vee a}\big) + \tfrac a{d_-\vee a}\sqrt{1-\big(\tfrac a{d_-\vee a}\big)^2}\Big)\\
	&\le C d_+\delta_d^{1/2}.
	\end{align*}	
\end{lemma}

\begin{proof}
	Recall that 
	\begin{equation*}
	D_{d_-,d_+}(x,y) = (B_{d_+}^+(x,y)\cup B_{d_+}^-(x,y)) \backslash (B_{d_-\vee a}^+(x,y)\cup B_{d_- \vee a}^-(x,y)).
	\end{equation*}
	The line through $x$ and $y$ cuts the disk $B_ d^+(x,y)$ into two parts. The area of the larger part is given by  
	\begin{align*}
	&d^2\Big(\pi -  \arcsin(\tfrac a d) + \tfrac a d\sqrt{1-(\tfrac a d)^2}\Big). 
	\end{align*}
	$D_{d_-,d_+}(x,y)$ is the union of two such sets of radius $d_+$ from which
	we remove two sets of the same type with radius $d_-\vee a$ from the interior. This yields the formula for the area.
	
	The inequality follows from
	\begin{align*}
	d_+^2-(d_-\vee a)^2{}&\le 2d_+ \delta_d, \\
	a(d_+-d_-\vee a)\sqrt{1-\big(\tfrac a{d_-\vee a}\big)^2}{}&\le d_+\delta_d,
	\end{align*}
	and, using \eqref{sqrtIneq} and \eqref{arcsinIneq},
	\begin{align*}
	 d_+^2{}&\Big(\arcsin\big(\tfrac a{d_-\vee a}\big) -  \arcsin\big(\tfrac a{d_+}\big) + \tfrac a{d_+}\Big(\sqrt{1-\big(\tfrac a{d_+}\big)^2}-\sqrt{1-\big(\tfrac a{d_-\vee a}\big)^2}\Big)\Big)\\
	 &\le
		d_+^2  \Big(C_0\sqrt{\tfrac a{d_-\vee a}-\tfrac a{d_+}}  + \tfrac a{d_+}\sqrt{(\tfrac a{d_-\vee a})^2 -(\tfrac a{d_+})^2}\Big) \\
	&\le C_1d_+^{3/2} 	\delta_d^{1/2}.
	\end{align*}
	
\end{proof}


\begin{lemma}\label{OneEdge}
	Let $0\le b_-<b_+ \le \tf$ and $0\le d_-<d_+ \le \tf$ and let $B_M$ be a disk of radius $M$. Then, there is a constant $C>0$ such that
	\begin{align*}
	\int_{B_M^3} \one_{(b_-,b_+]}\Big(\tfrac{|y_1-y_2|}2\Big) \one_{D_{d_-,d_+}(y_1,y_2)}(y_3) \d y_3 \d y_2 \d y_1 {}&\le C  |B_M| d_+^2 (\delta_b\vee \delta_d)^{\frac12}\delta_b\wedge \delta_d.
	\end{align*}
\end{lemma}

\begin{proof}

	Integration with respect to $y_3$ yields:
	\begin{align*}
	&\int_{B^3_M} \one_{(b_-,b_+]}\Big(\tfrac{|y_1-y_2|}2\Big) \one_{D_{d_-,d_+}(y_1,y_2)}(y_3) \d y_3 \d y_2 \d y_1\\
	& \qquad \le \int_{B^2_M} \one_{(b_-,b_+]}\Big(\tfrac{|y_1-y_2|}2\Big) |D_{d_-,d_+}(y_1,y_2)|  \d y_2 \d y_1.
	\end{align*}
	When $\delta_b\le \delta_d$, the claim follows directly from Lemma \ref{AreaDLem}. Otherwise, letting $a=|y_1-y_2|/2$, we split the integral in two terms according to whether $a<d_-$ or $a\geq d_-$. Applying Lemma \ref{AreaDLem} yields the bound
	\begin{align}\nonumber
	{}&C_1 |B_M| \bigg(\int_{b_-\wedge d_-}^{b_+\wedge d_-} a \Big(d_+^2\Big(\pi -  \arcsin\big(\tfrac a{d_+}\big) + \tfrac a{d_+}\sqrt{1-\big(\tfrac a{d_+}\big)^2}\Big)\\ \label{I1}
	& \qquad - 
	d_-^2\Big(\pi -  \arcsin\big(\tfrac a{d_-}\big) + \tfrac a{d_-}\sqrt{1-\big(\tfrac a{d_-}\big)^2}\Big) \Big)  \d a \\ 
	&+\int_{b_-\vee d_-}^{b_+\wedge d_+} a \Big(d_+^2\Big(\pi -  \arcsin(\tfrac a{d_+}) + \tfrac a{d_+}\sqrt{1-(\tfrac a{d_+})^2}\Big)  - a^2\tfrac{\pi}{2}\Big) \d a\bigg). \label{I2}
	\end{align}
	
	To bound \eqref{I1}, we apply the mean value theorem and perform the integration to obtain the bound 
	\begin{align*}
	C_1 |B_M| d_+^3 {}&  \int_{b_-\wedge d_-}^{b_+\wedge d_-}  \bigg( \tfrac1{\sqrt{1-\big(\tfrac a{d_-}\big)^2}}\big(\tfrac a{d_-} - \tfrac a{d_+}\big) + \tfrac{a^2}{d_+d_-}\tfrac1{\sqrt{1-\big(\tfrac a{d_-}\big)^2}}\big(\tfrac a{d_-} - \tfrac a{d_+}\big) \bigg)\d a\\
	&\le 2C_1 |B_M| d_+^2\delta_d \int_{b_-\wedge d_-}^{b_+\wedge d_-}   \tfrac a{\sqrt{d_-^2 -a^2}}\d a\\
	&=   2C_1 |B_M|d_+^2\delta_d \Big(\sqrt{d_-^2 -(b_-\wedge d_-)^2}-\sqrt{d_-^2 -(b_+\wedge d_-)^2}\Big)\\
	&\le   4C_1 |B_M|\tf d_+^2 \delta_d \delta_b^{1/2}.
	\end{align*}
	To bound \eqref{I2}, we bound the integrand using Lemma \ref{AreaDLem} and note that  
	\begin{equation*}
	|b_+\wedge d_+- b_-\vee d_-| \le \delta_d\wedge \delta_b.
	\end{equation*}
	This proves the claim when $\delta_d\le \delta_b$.
\end{proof}


\begin{lemma}\label{specialconfig}
	 Let $B_M$ be a disk of radius $M > 0$. There is a constant $C>0$ such that for all $b_-,b_+,b_-',b_+',d_-,d_+ \in [0,\tf]$ with $d_-<d_+$ and either $b_-<b_+=b_-'<b_+'$ or $b_-=b_-'$ and $b_+=b_+'$,
	\begin{align*}
	&\int_{B^4_M} \one_{ (b_-, b_+]\times  (b_-', b_+']\times (0,b_+\vee b_+'] }(\tfrac{|x_1 - x_2|}2,\tfrac{|x_1 - x_3|}2,\tfrac{|x_2 - x_3|}2)\\
	&\qquad \times \one_{ D_{d_-, d_+}( x_2,x_3)}(y_1)\d y_1 \d x_1 \d x_2 \d x_3
	\\ 
	&\quad \le C|B_M|\delta_b\wedge \delta_{b'} (\delta_b \vee \delta_b')^{3/4} \delta_d^{3/4}.
	\end{align*}
\end{lemma}

\begin{proof}

	We may assume $d_+ > 3\delta_d\vee8\sqrt{\tf(\delta_b + \delta_{b'})}$. Indeed, if $d_+\le 3\delta_d$, we can show the claim by first integrating  with respect to $y_1$, then using that by Lemma \ref{OneEdge},
	\begin{equation*}
	|D_{d_-, d_+}(x_2,x_3)|\le C_1d_+^2 \le 9C_1 \delta_d^2,
	\end{equation*}
	and finally integrating with respect to $x_2$ and $x_3$ to provide a factor $|B_M|\delta_b\delta_{b'}$. 
	If $d_+\le 8\sqrt{\tf(\delta_b + \delta_{b'})}$, we first integrate with respect to $x_1$, which yields the area of $A_{2b_-,2b_+}(x_2)\cap A_{2b_-',2b_+'}(x_3)$. This is bounded by  $C_2\delta_{b}\wedge \delta_{b'}$, and by Lemma~\ref{OneEdge} the remaining integral is bounded by 
	\begin{equation*}
	C_3|B_M|d_+^2\delta_d \le 64C_3|B_M|\tf (\delta_b + \delta_{b'})\delta_d \le  128 C_3|B_M|\tf (\delta_b\vee \delta_{b'}) \delta_d.
	\end{equation*}
	
	Let $a=|x_2 - x_3|/2$. We write the integral as a sum of three terms corresponding to whether 
	I: $a<d_+/4$, 
	II: $d_+/4 \le a < b_-\wedge b_-'$, or 
	III: $b_-\wedge b_-' \le a\le b_+\vee b_+'$. 
	
	Term I: We first integrate with respect to $y_1$. Since 
	\begin{equation*}
	\tfrac a{d_-}\le \tfrac{d_+}{4d_-} = \tfrac{d_- + \delta_d}{4d_-} \le \tfrac34,
	\end{equation*}
	the mean value theorem applied to the formula in Lemma \ref{AreaDLem} implies that $|D_{d_-, d_+}( x_2,x_3)|\le C_4\delta_d$. We then integrate with respect to $x_2$ and $x_3$ to obtain the bound $C_5|B_M|\delta_b\delta_{b'}\delta_d$.
	
	Term II:
	When $d_+/4 \le a \le b_-\wedge b_-'$, we first integrate with respect to $x_1$ to obtain the area of $ A_{2b_-,2b_+}(x_2)\cap A_{2b_-',2b_+'}(x_3)$. To bound term II, we need to explicitly compute this area. For this, we first compute the area $A_a(b_1,b_2)$ of the intersection $B_{2b_1}(x_2)\cap B_{2b_2}(x_3)$ where $b_1,b_2 \in \{b_+,b_-,b_+',b_-'\}$.
	By the assumption on $d_+$, 
	\begin{equation}\label{asquared}
	a^2\geq d_+^2/16 \geq 4 \tf (\delta_b+\delta_{b'})\geq 2(b_1^2-b_2^2).
	\end{equation}
	This ensures that the line containing the two points where the boundaries of the disks $B_{2b_1}(x_2)$ and $B_{2b_2}(x_3)$ meet separates $x_2$ and $x_3$. The area of $B_{2b_1}(x_2)\cap B_{2b_2}(x_3)$ is 
	\begin{align*}
	\mathcal{A}_a(b_1,b_2)={}&4\bigg(b_1^2\arccos\Big(\tfrac{a^2+b_1^2-b_2^2}{2ab_1}\Big)+b_2^2\arccos\Big(\tfrac{a^2+b_2^2-b_1^2}{2ab_2}\Big)\\
	& - b_1\tfrac{a^2+b_1^2-b_2^2}{2a}\Big(1-\Big(\tfrac{a^2+b_1^2-b_2^2}{2ab_1}\Big)^2\Big)^{1/2} \\
	&- b_2\tfrac{a^2+b_2^2-b_1^2}{2a}\Big(1-\Big(\tfrac{a^2+b_2^2-b_1^2}{2ab_2}\Big)^2\Big)^{1/2}\bigg).
	\end{align*}
	The area of $ A_{2b_-,2b_+}(x_2)\cap A_{2b_-',2b_+'}(x_3)$ is given by
	\begin{align}\nonumber
	\mathcal{A}_a(b_+,b_+'){}& + \mathcal{A}_a(b_-,b_-') - \mathcal{A}_a(b_+,b_-') - \mathcal{A}_a(b_+',b_-) \\
	{}&= \int_{b_-}^{b_+} \int_{b_-'}^{b_+'}\tfrac{\partial^2}{\partial b_1 \partial b_2} \mathcal{A}_a(b_1,b_2) \d b_1 \d b_2.\label{Ab1b2} 
	\end{align}
	It is a straightforward computation to see that $\tfrac{\partial^2}{\partial b_1 \partial b_2} \mathcal{A}_a(b_1,b_2)$  is uniformly bounded by $C_6/d_+^2$ on the set of $a,b_1,b_2 \le \tf$ satisfying \eqref{asquared} and $d_+/4\le a \le b_1\wedge b_2$. In particular, \eqref{asquared} guarantees that 
	\begin{equation*}
	\frac{a^2+b_1^2-b_2^2}{2ab_1} \le \frac{3a}{4b_1} \le \frac34,
	\end{equation*}
	such that $\arccos$ and $x\mapsto \sqrt{1-x^2}$ have bounded derivatives for the relevant values of $x$. It follows that  \eqref{Ab1b2} is bounded by $C_7\delta_b \delta_{b'}/d_+^2$.  The remaining integral is of order $|B_M|d_+^2 \delta_d$ by Lemma \ref{OneEdge}, which yields the  appropriate bound.
	
	Term III: In this case, we first integrate with respect to $x_1$ providing a factor $\delta_b\wedge \delta_{b'}$. The remaining integral is bounded using Lemma \ref{OneEdge}.
\end{proof}

The fourth lemma allows us to analyze which point configurations can cause the birth and death of $M$-bounded features. To state it, we recall the $\alpha$-complex associated with a locally finite point set $\X\subseteq \R^2$, see e.g.\ \cite[Sec.~III.4]{edHar} for details. It is built from the Delaunay triangulation, which is a triangulation of the plane with vertex set $\X$. For $r>0$, $\alpha_r(\X )$ is the union of all edges in the Delaunay triangulation with length at most $2r$ and all triangles such that the three balls of radius $r$ centered at its vertices cover the triangle. Then $\alpha_r(\X)\subseteq U_r(\X)$ and the inclusion is a homotopy equivalence, i.e.\ it preserves the topology.

\begin{lemma}\label{alpha}
	Let $\X \subseteq \R^2$ be locally finite.		
	\begin{itemize}
	\item[(i)]  Each connected component of $\R^2 \backslash \alpha_r(\X)$ contains at most one $M$-bounded connected component of $\R^2\backslash U_r(\X)$. 
	\item[(ii)] If an $M$-bounded loop is born at time $b$ because two balls centered at $x_1,x_2$ meet, then there is an edge of length $2b$ joining $x_1,x_2$ in the $\alpha$-complex. 
	\item[(iii)] If an $M$-bounded feature dies at time $d$ because exactly three balls centered at points $y_1,y_2,y_3$ meet, then $y_1,y_2,y_3$ form a triangle with no obtuse angle in the $\alpha$-complex. 
	\end{itemize}
\end{lemma}

\begin{proof}
	The analogous statements hold for unbounded loops by the homotopy equivalence between the $\alpha$-complex and the union of balls. (i) follows because any $M$-bounded loop is also an unbounded loop.  An $M$-bounded feature is either born the same way as the corresponding unbounded component or when two balls meet to split off a component. In both cases, some unbounded loop is born by the merging, and hence an edge is added to the $\alpha$-complex. This shows (ii). When an $M$-bounded loop dies, so does the corresponding unbounded loop, hence (iii) is clear.	
\end{proof}


We are now ready to prove Lemma \ref{geom1Lem}.

\begin{proof}[Proof of Lemma \ref{geom1Lem}]

	{\bf Proof of \eqref{moment1}.}
	Stationarity and Equation \eqref{redPalmEq} yield
	\begin{align}
	\begin{split}\label{disint}
	&\int_{B^p} {\P}_{o, \zz}(E_o) \rho^{(p + 1)}(o, \zz) \d \zz \\
	&\quad = \int_{[0, 1]^2}\int_{(B + x)^p} {\P}_{x, \zz}(E_x) \rho^{(p + 1)}(x, \zz) \d \zz \d x\\
		&\quad = \E \bigg[ \sum_{(x, \zz)\in \PP^{p + 1}_{ \ne}} \one_{[0, 1]^2}(x) \one_{(B + x)^p}(\zz) \one_{E_x}\bigg] .
	\end{split}
	\end{align}
	In the following, we let $\yy = (y_1, y_2, y_3)$, and 
	\begin{equation*}
	g(x_1, x_2, \yy) = \one_{ (b_-, b_+]}\Big(\tfrac{|x_1 - x_2|}2\Big) \one_{ D_{d_-, d_+}(y_1, y_2)}(y_3)
	\end{equation*}
	for simplicity. By definition of $E_{x} $, \eqref{disint} is bounded by
	\begin{align}
	\begin{split}\label{disjointSum}
	&\E \bigg[ \sum_{x_1\in \PP}\PP(B + x_1)^p \one_{[0, 1]^2}(x_1) \sum_{x_2\in \PP} \sum_{\yy \in \PP ^3_{\ne} } \one_{B_M(x_1)^3}(\yy) g(x_1, x_2, \yy) \bigg] \\	
	&\quad = \E \bigg[\sum_{(x_1, x_2, \yy) \in \PP^5_{\ne} } \PP(B + x_1)^p \one_{[0, 1]^2}(x_1) \one_{B_M(x_1)^3}(\yy) g(x_1, x_2, \yy)\bigg] \\
	&\quad \quad + 3 \E\bigg[ \sum_{(x_1, \yy) \in \PP^4_{\ne} } \PP(B + y_1)^p \one_{[0, 1]^2}(y_1)\one_{B_M(y_1)^2}(y_2, y_3) g(x_1, y_1, \yy) \bigg]\\
	&\quad \quad + 3\E \bigg[\sum_{(x_1, \yy) \in \PP^4_{\ne} } \PP(B + x_1)^p \one_{[0, 1]^2}(x_1) \one_{B_M(x_1)^3}(\yy) g(x_1, y_1, \yy) \bigg]\\
	&\quad \quad + 6 \E \bigg[\sum_{\yy \in \PP^3_{\ne} } \PP(B + y_2)^p\one_{[0, 1]^2}(y_2)\one_{B_M(y_2)^2}(y_1, y_3) g(y_1, y_2, \yy)\bigg].
	\end{split}
	\end{align}
	Here, we have used that $g(x_1, x_2, \yy)$ is symmetric in $x_1$ and $x_2$ and in $y_1$, $y_2$, and $y_3$. Applying \eqref{redPalmEq} again, we may bound the last term in \eqref{disjointSum} by
	\begin{align}
	\begin{split}\label{3integrals}
	6 \int_{B_{M + 2}^3} \E_{\yy}[\PP(B + y_1)^p] g(y_1, y_2, \yy) \rho^{(3)}(\yy) \d \yy,
	\end{split}
	\end{align}
	since $b_+ \le M$. The remaining terms are treated similarly. Now choose a covering $B + x_1\subseteq \bigcup_{i \le \ell} W_1^{(i)}$, where each $W_1^{(i)}$ is a translation of $W_1$ and such that $\ell\le C_1|B|$ for some $C_1$ independent of $K$ (for instance using that $B_K \subseteq W_{4\lceil K \rceil^2}$). Then, by the moment condition {\bf (M)} for $\xx = (x_1, \dots, x_k)$, 
	\begin{align*}
	\E_{\xx}\big[\PP(B_K+x_1)^p \big] {}& \le \sup_{\xx \in \R^{2k}}\E_{\xx}\Big[\PP\Big(\bigcup_{i = 1}^\ell W_1^{(i)} \Big)^p \Big]\\
	&\quad \le \ell^p\sum_{i \le \ell} \sup_{\xx \in \R^{2k}}\E_{\xx}\big[\PP\big( W_1^{(i)}\big)^p\big]\\
	&\quad \le \ell^p\sum_{i \le \ell}\sup_{\xx \in \R^{2k}} \E_{\xx}^!\big[\big(\PP\big( W_1^{(i)}\big) + k\big)^{p\vee k}\big]\\
	&\quad \le C_2 \ell^{p + 1} \Big(\sup_{\xx \in \R^{2k}} \E_{\xx}^!\big[\PP( W_1)^{p\vee k}\big] + k^{p\vee k}\Big)\\
	&\quad \le C_3 |B|^{p + 1}.
	\end{align*}
	We apply this in \eqref{3integrals} together with Lemma \ref{OneEdge}. Since each $\rho^{(k)}$ is bounded according to the assumption of fast decay of correlations, we obtain the bound $C_4 |B|^{p + 1} |E|^{1/2+\e}$.

	{\bf Proof of \eqref{moment2}. }
	In the following, we use the notation
		\begin{equation*}
	g'(x_1, x_2, \yy) = \one_{ (b'_-, b'_+] }\Big(\tfrac{|x_1 - x_2|}2\Big) \one_{ D_{d'_-, d'_+}(y_1, y_2)}(y_3).
	\end{equation*}
		Note that since the blocks $E$ and $E'$ are neighboring, the features in $E$ and $E'$ are different.
Putting $\xx = (x_1, x_2, x_3)$, we now expand as in \eqref{disint}
	\begin{align}
	\begin{split}\label{9pointsum}
		&\int_{B^p} \P_{o, \zz}(E_{o, o}'') \rho^{(p + 1)}(o, \zz) \d \zz\\
		&\quad \le \E \bigg[\sum_{\substack{x_1\in \PP\cap [0, 1]^2\\ (x_2, x_3)\in \PP^2_{\ne} }}  \sum_{\substack{\yy, \yy' \in \PP_{\ne}^3 \cap B_M(x_1)^3 \\ \yy \ne \yy'}} \hspace{-.5cm}\PP(B + x_1)^p g(x_1, x_2, \yy) g'(x_1, x_3, \yy')\one_A(\xx,\yy,\yy')\bigg]\\
	&\quad \le \E \bigg[ \sum_{\substack{\xx, \yy, \yy'\in \PP^3_{\ne} \cap B_{M + 2}^3\\ \yy\ne \yy'}} \PP(B + x_1)^p
	g(x_1, x_2, \yy) g'(x_1, x_3, \yy')\one_A(\xx,\yy,\yy')\bigg].
	\end{split}
	\end{align}
		The condition $x_2\ne x_3$ comes from the fact that $x_1$ can give birth to at most one feature when connecting to another point, and since $E$ and $E'$ are neighboring, $x_2$ and $x_3$ correspond to different features. Similarly, $\yy'\ne\yy$ comes from the fact that a triangle can kill at most one feature.
	
	The event $A$ excludes certain point configurations that are not possible. If the triangles formed by $\yy$ and $\yy'$ share an edge, and the vertices of this edge coincide with $x_2$ and $x_3$, then $|x_2 - x_3|>2(b_+\vee b_+')$ is not allowed. Indeed, it follows from  Lemma \ref{alpha} that the triangles correspond to the same feature in the $\alpha$-complex until $x_2$ and $x_3$ are joined. Thus, this must happen before both triangles are born, that is, at the latest at time $b_+\vee b_+'$. Moreover, if the two triangles share an edge, then the two points in $\yy,\yy'$ not lying on this edge cannot be equal to $x_1$ and $x_2$ or to $x_1$ and $x_3$, as this would lead to crossing edges in the $\alpha$-complex by Lemma \ref{alpha} (since the triangles formed by $\yy,\yy'$ cannot have any obtuse angles).

We now write the sum in \eqref{9pointsum} as a sum where each term is a sum over $\PP_{\ne}^k$, $4\le k\le 9$, as in \eqref{disjointSum}. Each such term comes from grouping $\xx, \yy, \yy'$ into sets of equal points. Consider for illustration the term corresponding to the situation $x_2 = y_1', x_1 = y_2 = y_2', x_3= y_3 = y_3'$. The sum is handled as in the proof of Lemma \ref{geom1Lem}\eqref{moment1} by applying \eqref{redPalmEq} and bounding the involved Palm means. For this special point configuration, it is sufficient to bound $\one_A$ by 1.	
\begin{align*}
	&\frac1{|B|^{p + 1}}\E\Big[ \sum_{(\xx, y_1)\in \PP^4_{\ne} \cap B_{M + 2}^4} \PP(B + x_1)^p g(x_1, x_2, y_1, x_1, x_3) g'(x_1, x_2, \xx)
	\Big]\\	
	&\quad \le C_5\int_{B_{M + 2}^4} g(x_1, x_2, y_1, x_1, x_3) g'(x_1, x_3, \xx) \d y_1 \d \xx.
\end{align*}
Now, we apply the H\"older inequality with $\frac1{q_1} + \frac1{q_2} = 1$ to obtain the bound
	\begin{align}
	\begin{split} \label{holderStep}
	& C_5 \Big[\int_{B_{M + 2}^4}\hspace{-.4cm} \one_{ (b_-, b_+] \times (b_-', b_+']}\Big(\tfrac{|x_1 - x_2|}2, \tfrac{|x_1 - x_3|}2\Big) \one_{ D_{d_-, d_+}(x_1, x_3)}(y_1)\d y_1\d \xx\Big]^{\frac1{q_1}}\\
		&\quad \times \Big[\int_{B_{M + 2}^4}\hspace{-.4cm} \one_{ D_{d_-, d_+}(x_1, x_3)}(y_1)\one_{ (b_-', b_+']}\Big(\tfrac{|x_1 - x_3|}2\Big) \one_{ D_{d_-', d_+'}(x_1, x_2)}(x_3)\d y_1\d \xx\Big]^{\frac1{q_2}}.
	\end{split}
	\end{align}
 	In the first integral, we first integrate with respect to $x_2$ and then apply Lemma \ref{OneEdge}, while in the second integral we first integrate with respect to $y_1$ and use the bound in Lemma \ref{AreaDLem} and then apply Lemma \ref{OneEdge} again. 	Next we use that $E$ and $E'$ are neighboring blocks so that either $\delta_b=\delta_{b'} $ or $\delta_ d=\delta_{d'}$.
	
	When $\delta_b=\delta_{b'}$, we get the bound
	\begin{equation}\label{boundb}
	C_6 \big(\delta_{b}(\delta_{b'}\delta_d)^{\frac34}\big)^{\frac1{q_1}}
	\big(\delta_d^{\frac12}(\delta_{b'}\delta_{d'})^{\frac34}\big)^{\frac1{q_2}}  = C_6\delta_{b}^{\frac34 + \frac1{q_1}}\delta_d^{\frac34 \cdot \frac1{q_1} + \frac12 \cdot \frac1{q_2}}\delta_{d'}^{\frac34 \cdot \frac1{q_2}}, 
	\end{equation}
	so we take $1/q_1 > 1/4$ and $1/q_2 > 2/3$.
	
	When $\delta_d=\delta_{d'}$, we use Lemma \ref{OneEdge} to get the bound
	\begin{equation}\label{boundd}
	C_7 \big(\delta_{b}(\delta_{b'}\delta_d)^{\frac34}\big)^{\frac1{q_1}}
	\big(\delta_d^{\frac12}\delta_{b'}^{\frac12}\delta_{d'}\big)^{\frac1{q_2}}  = C_7\delta_{b}^{ \frac1{q_1}}\delta_{b'}^{\frac34 \cdot \frac1{q_1} + \frac12 \cdot \frac1{q_2}}\delta_ d^{\frac34 \cdot \frac1{q_1} + \frac32 \cdot \frac1{q_2}}, 
	\end{equation}
	so we take $1/q_1 > 1/2$ and $1/q_2 > 1/3$.
	
For a general term, note that there are at least four different points among $\yy,\yy'$, so one of them, say $y_1$, cannot be equal to any of $\xx$. We consider two cases: 
\begin{itemize}
	\item[I] $y_1$ is not among $y_1',y_2',y_3'$.
		\item[II] $y_1=y_1'$, $y_2=y_2'$, and $y_3=x_2$ and $y_3'=x_3$.
\end{itemize}
Since we no longer keep track of which edge kills which triangle, all possible point configurations allowed by $A$ fall into one of the above cases after possibly renaming the variables.

 In particular, if $y_1=y_1'$ and the points $y_2,y_3,y_2',y_3'$ are all different, one of them cannot be any of $x_1,x_2,x_3$, and we could have taken this as $y_1$ and be in Case I.  If $y_1=y_1'$, $y_2=y_2'$ and, say, $y_3$ is not any of $x_1$, $x_2$, $x_3$, we could have chosen $y_3$ as $y_1$ and be in Case I.

We further divide the Case I configurations allowed by $A$ into the following two sub-cases that have to be treated separately:
\begin{itemize}
	\item[Ia.]  $x_3$ is not any of $y_2,y_3$.
	\item[Ib.]  $x_2=y_2=y_2'$, $x_3=y_3=y_3'$, $|x_2 - x_3|/2\le b_+\vee b_+'$.
\end{itemize}
Again, after renaming the variables, we are always in one of the two sub-cases. 

Case Ia:  We apply the H\"older inequality to
\begin{align}\nonumber
&\one_{ (b_-, b_+]}\Big(\tfrac{|x_1 - x_2|}2\Big) \one_{ D_{d_-, d_+}(y_1, y_2)}(y_3)\one_{ (b_-', b_+']}\Big(\tfrac{|x_1 - x_3|}2\Big) \one_{ D_{d_-', d_+'}(y_1', y_2')}(y_3')\\
&\quad =  \one_{ (b_-, b_+]}\Big(\tfrac{|x_1 - x_2|}2\Big)\one_{ (b_-', b_+']}\Big(\tfrac{|x_1 - x_3|}2\Big) \one_{ D_{d_-, d_+}(y_1, y_2)}(y_3)\label{indic1}\\
&\quad \quad\times\one_{ D_{d_-, d_+}(y_1, y_2)}(y_3)\one_{ (b_-', b_+']}\Big(\tfrac{|x_1 - x_3|}2\Big) \one_{ D_{d_-', d_+'}(y_1', y_2')}(y_3'). \nonumber
\end{align}
The first factor is integrated with respect to $x_3$ and the remaining integral is bounded using Lemma \ref{OneEdge}.  The second factor is first integrated wrt. $y_1$, the result is bounded using Lemma \ref{AreaDLem}, and the remaining integral is bounded using Lemma \ref{OneEdge}. The rest of the argument proceeds as in the special case treated above.

Case Ib: The claim follows by applying the Hölder inequality to \eqref{indic1} and arguing as in Case Ia using Lemma \ref{specialconfig} to bound the first integral.

Case II:  
We apply the H\"{o}lder inequality exactly as in \eqref{indic1} and argue as in Case Ia, except that the second integral is first integrated with respect to $y_3$ rather than $y_1$.

	{\bf Proof of \eqref{moment3}.} As in \eqref{disint}, we find
	\begin{align*}
		&\int_{B^{p + 1}} \P_{o, z', \zz}(E_{o,z'}'') \rho^{(p+2)}(o, z', \zz) \d z' \d \zz\\
		&\quad \le \E \bigg[ \sum_{(x, z')\in \PP_{n \ne}^2} \PP(B + x)^p \one_{[0, 1]^2}(x)\one_{B + x}(z') \one_{E_{x, z'}''}\bigg] \\
		&\quad \le \E \bigg[\sum_{(x_1, z')\in \PP_{\ne}^2 } \sum_{(x_2, x_2')\in \PP^2}  \sum_{\yy \in \PP_{\ne}^3 \cap B_M(x_1)^3} \sum_{\substack{ \yy'\in \PP_{\ne}^3 \cap B_M(z')^3\\ \yy \ne \yy'}}  \PP(B + x_1)^p \one_{[0, 1]^2}(x_1)\\
	&\quad\quad \times \one_{B + x_1}(z') g(x_1, x_2, \yy) g'(z', x_2', \yy')\one_{\tilde{A}}(x_1,z',x_2,x_2',\yy,\yy')\bigg].
	\end{align*}
		The set $\tilde{A}$ consists of tuples of points $ (x_1,x_2,x_3,x_4,\yy,\yy')\in \R^{20}$ and, similar to $A$, it excludes certain configurations of the points  $ (x_1,x_2,x_3,x_4,\yy,\yy')$ that are not allowed by Lemma \ref{alpha}. If the triangles formed by $\yy$ and $\yy'$ share an edge,  then the length of this edge must be at most $2(b_+\vee b_+')$. Moreover, if the two triangles share an edge, then the two points in $\yy,\yy'$ not lying on this edge cannot be equal to $x_1$ and $x_3$ or to $x_2$ and $x_4$.
		 
	The contribution from the cases where two of the points $x_1, x_2, x_2', z'$ are identical is bounded by
	\begin{align*}
	\begin{split}
	&\E \bigg[  \sum_{\substack{\xx, \yy, \yy'\in \PP_{\ne}^3\cap B_{2M + 2}^3\\ \yy'\ne\yy }}  \PP(B + x_1)^p g(x_1, x_2, \yy) g'(x_1, x_3, \yy')\one_{\tilde{A}}(x_1,x_2,x_1,x_3,\yy,\yy') \bigg], 
	\end{split}
	\end{align*}
	which is handled exactly as in the proof of Lemma  \ref{geom1Lem}\eqref{moment2}. Thus, it remains to treat the terms where $x_1, x_2, x_2', z'$ are all different. Therefore, if we put $\xx = (x_1, x_2, x_3, x_4)$, we must bound
	\begin{align*}
		&\E \bigg[\sum_{\substack{\xx \in \PP_{\ne}^4 \\ x_1 \in [0, 1]^2}} \sum_{\yy \in \PP_{\ne}^3 \cap B_M(x_1)^3} \sum_{\substack{ \yy'\in \PP_{\ne}^3 \cap B_M(x_2)^3\\ \yy \ne \yy'}}  \PP(B + x_1)^p  \one_{B + x_1}(x_2)\\
		&\quad \times g(x_1, x_3, \yy)g'(x_2, x_4, \yy')\one_{\tilde{A}}(\xx,\yy,\yy') \bigg].
	\end{align*}
	The rest of the proof proceeds as the proof of  Lemma  \ref{geom1Lem}\eqref{moment2} by suitable applications of the H\"{o}lder inequality. 
We divide into two cases according to whether all points in $\yy$, $\yy'$ are one of $\xx$ or not. After renaming the variables, we may assume
\begin{itemize}
	\item[I] $y_1=y_1'=x_1$, $y_2=y_2'=x_2$, $y_3=x_3$, and $y_3'=x_4$, or
	\item[II] $y_1$ is not any of $\xx$.
\end{itemize}
Notice that in Case I we exclude the case $y_1=y_1'=x_1$, $y_2=y_2'=x_3$, $y_3=x_2$, and $y_3'=x_4$ because it was excluded by definition of $\tilde{A}$.
After renaming variables, Case II is divided into
\begin{itemize}
	\item[IIa] $y_1$ is not any of $\xx$ or $\yy'$, and $x_1$ is not any of $y_2,y_3$.
	\item[IIb] $y_1=y_1'$ and $y_1$ is not any of $\xx$, $y_2=y_2'\ne x_3$, $y_3=x_1$.
	\item[IIc] $y_1=y_1'$, $y_2=x_2$, $y_3=x_4$, $y_2'=x_1$, $y_3'=x_3$.
\item[IId] $y_1=y_1'$, $y_2=x_1$, $y_3=x_2$, $y_2'=x_3$, $y_3'=x_4$.
\end{itemize}
In Case IIa, $y_1$ is not one of $\yy'$, while in Case IIb, IIc, and IId it is. 
 Case IIb corresponds to the situation in which the triangles formed by $\yy,\yy'$ share an edge, while in Case IIc and IId  they share only one vertex. In Case IIc, each triangle contains one of the edges joining $x_1$ to $x_3$ and $x_2$ to $x_4$, while in Case IId they do not. 

Case I: When $\delta_b=\delta_{b'}$, we first write
\begin{align}\nonumber
&\one_{ (b_-, b_+]}\Big(\tfrac{|x_1 - x_3|}2\Big) \one_{ D_{d_-, d_+}(y_1, y_2)}(y_3)\one_{ (b_-', b_+']}\Big(\tfrac{|x_2-x_4|}2\Big) \one_{ D_{d_-', d_+'}(y_1', y_2')}(y_3')\\ \label{2b1d}
&\quad = \one_{ (b_-, b_+]}\Big(\tfrac{|x_1 - x_3|}2\Big) \one_{ D_{d_-, d_+}(y_1, y_2)}(y_3) \one_{ (b_-', b_+']}\Big(\tfrac{|x_2-x_4|}2\Big)\\ \label{1b2d}
&\quad \quad\times \one_{ (b_-', b_+']}\Big(\tfrac{|x_2 - x_4|}2\Big) \one_{ D_{d_-, d_+}(y_1, y_2)}(y_3)\one_{ D_{d_-', d_+'}(y_1', y_2')}(y_3').
\end{align}
We then apply the H\"{o}lder inequality.
Integrating first with respect to $x_4$ and then $y_2$ in \eqref{2b1d} and integrating with respect to $y_3$ first in \eqref{1b2d} yields a bound of order
\begin{equation*}
(\delta_{b'} (\delta_{b}\delta_d)^{\frac34})^{\frac1{q_1}}(\delta_d (\delta_{b'}\delta_{d'})^{\frac34})^{\frac1{q_2}}.
\end{equation*}
This is the same as \eqref{boundb} since $\delta_b=\delta_{b'}$. When $\delta_d=\delta_{d'}$, we replace $\one_{ D_{d_-, d_+}(y_1, y_2)}(y_3)$ by $\one_{ D_{d_-', d_+'}(y_1', y_2')}(y_3')$ in \eqref{2b1d}, to obtain a bound of order 
\begin{equation*}
(\delta_{b} (\delta_{b'}\delta_{d'})^{\frac34})^{\frac1{q_1}}(\delta_d^{\frac12} \delta_{b'}^{\frac12}\delta_{d'})^{\frac1{q_2}},
\end{equation*}
 which reduces to the same form as  \eqref{boundd}. 

Case IIa: We apply the H\"{o}lder inequality to \eqref{2b1d}--\eqref{1b2d} and integrate first with respect to $x_1$ and then $y_1$ in \eqref{2b1d} and with respect to $y_1$ first in \eqref{1b2d}. The remaining argument proceeds as in the proof of Lemma \ref{geom1Lem}\eqref{moment2} Ia.

Case IIb:  We apply the H\"{o}lder inequality to \eqref{2b1d}--\eqref{1b2d} and integrate first with respect to $x_3$ and then $y_1$ in \eqref{2b1d} and with respect to $y_3$ first in \eqref{1b2d} and argue as in the proof of Lemma \ref{geom1Lem}\eqref{moment2} Ia.

Case IIc: In \eqref{2b1d}, we first integrate with respect to $x_1$. In \eqref{1b2d}, we first integrate with respect to $y_2'$ and $y_3'$ to obtain a factor $\delta_{d'}$. Then we integrate with respect to $y_1$ and $x_2$ and apply Lemma \ref{OneEdge} to obtain a factor $\delta_d^{1/2}\delta_{b'}$. The resulting bounds are stricter than \eqref{boundb} and \eqref{boundd}.

Case IId: Here we integrate \eqref{2b1d} with respect to $x_3$ first and then $y_1$ while \eqref{1b2d} is integrated first with respect to $y_2$ and then $y_1$. 

In all cases treated above, a minor difference to \eqref{holderStep} is that the integration domains are slightly more complicated due to the indicator $\one_{B + x_1}(x_2)$. However, it contributes at most a factor $C_7 |B|$ to the bound, and this cancels when we divide by $|B|^{p+2}$.

\end{proof}

\end{document}